\documentclass[smallextended]{svjour3_CorrectedPostsubmission}
\usepackage{amsfonts}
\usepackage{graphicx}
\usepackage{epstopdf}
\usepackage{amsmath, amssymb, bm}
\usepackage[fixamsmath,disallowspaces]{mathtools}

\RequirePackage{fix-cm}
\smartqed  
\journalname{Multibody System Dynamics}
\begin{document}

\vspace{-10ex}\title{Screw and Lie Group Theory in Multibody Kinematics}
\subtitle{Motion Representation and Recursive Kinematics of Tree-Topology Systems}
\author{Andreas M\"uller}
\vspace{-10ex}

\date{
\,
}

\institute{A. M\"uller \at Institute of Robotics, Johannes Kepler University, Altenberger Str. 69, 4040 Linz, Austria, a.mueller@jku.at
}
\vspace{-10ex}

\titlerunning{Screw Theory in MBS Kinematics}
\maketitle

\begin{abstract}
After three decades of computational multibody system (MBS) dynamics,
current research is centered at the development of compact and user friendly
yet computationally efficient formulations for the analysis of complex MBS.
The key to this is a holistic geometric approach to the kinematics modeling
observing that the general motion of rigid bodies as well as the relative
motion due to technical joints are screw motions. Moreover, screw theory
provides the geometric setting and Lie group theory the analytic foundation
for an intuitive and compact MBS modeling. The inherent frame invariance of
this modeling approach gives rise to very efficient recursive $O\left(
n\right) $ algorithms, for which the so-called 'spatial operator algebra' is
one example, and allows for use of readily available geometric data. In this
paper three variants for describing the configuration of tree-topology MBS
in terms of relative coordinates, i.e. joint variables, are presented: the
standard formulation using body-fixed joint frames, a formulation without
joint frames, and a formulation without either joint or body-fixed reference
frames. This allows for describing the MBS kinematics without introducing
joint reference frames and therewith rendering the use of restrictive modeling convention, 
such as Denavit-Hartenberg parameters, redundant. Four
different definitions of twists are recalled and the corresponding
recursive expressions are derived. The
corresponding Jacobians and their factorization are derived. The aim of this
paper is to motivate the use of Lie group modeling and to provide a review
of the different formulations for the kinematics of tree-topology MBS in
terms of relative (joint) coordinates from the unifying perspective of screw
and Lie group theory.
\end{abstract}


\section{Introduction}

Computational multibody system (MBS) dynamics aims at mathematical
formulations and efficient computational algorithms for the kinetic analysis
of complex mechanical systems. At the same time the modeling process is
supposed to be intuitive and user friendly. Moreover, the efficiency of MBS
algorithms as well as the complexity of the actual modeling process is
largely determined by the way in which the kinematics is described. This
concerns the core issues of representing rigid body motions and of
describing the kinematics of technical joints. Both issues can be addressed
with concepts of screw and Lie group theory.

Spatial MBS perform complicated motions, and in general rigid bodies perform
screw motions that form a Lie group. Although the theory of screw motions is
well understood, screw theory has almost completely been ignored for MBS
modeling with only a few exceptions. The latter can be grouped into two
classes. The first class makes use of the fact that the velocity of a rigid
body is a screw, referred to as the \emph{twist}. The propagation of twists
within an MBS is thus described as a frame transformation of screw
coordinates. This gave rise to the so-called 'spatial vector' formulation
introduced in \cite{Featherstone1983,Featherstone2008}, and to the so-called
'spatial operator algebra' that was formalized in \cite{Rodriguez1991} and
used for $O\left( n\right) $ forward dynamics algorithms e.g. in \cite%
{Fijany1995,Jian1991,Jian1995,LillyOrin1991,Rodriguez1987,Rodriguez1992}.
Screw notations are also used in the formulations presented in \cite%
{Angeles2003,Legnani1,Legnani2,Uicker2013}. Further MBS formulations were
reported that use screw notations uncommon for the MBS community \cite%
{Gallardo2003}. All these approaches only exploit the algebraic properties
of screws as far as relevant for a compact handling of velocities,
accelerations, wrenches, and inertia. The second class goes one step further
by recognizing that finite motions form the Lie group $SE\left( 3\right) $
with the screw algebra as its Lie algebra $se\left( 3\right) $. Moreover,
screw theory provides the geometric setting and Lie group theory the
analytic foundation for an intuitive and efficient modeling of rigid body
mechanisms. Some of the first publications reporting Lie group formulations
of the kinematics of an open kinematic chain are \cite%
{Brockett1984,Hervet1978,Hervet1982} and \cite{Chevalier1991,Chevalier1994}.
In this context the term \emph{product of exponentials} (POE) is being used
since Brockett used it in \cite{Brockett1984}. Unfortunately these
publications have not reached the MBS community, presumably because of the
used mathematical concepts that differ from classical MBS formalisms. The
basic concept is the exponential mapping that determines the finite relative
motion of two adjacent bodies connected by a lower pair joint in terms of a
screw associated with the joint. The product of the exponential mappings for
all consecutive joints determines the overall motion of the chain. Within
this formulation twists are naturally represented as screws, and joint
motions are described in terms of screw coordinates.

Motivated by \cite{Brockett1984} Lie group formulations for MBS dynamics
were reported in a few publications, e.g. \cite%
{Mladenova2006,MUBOLieGroup,Park1994,ParkBobrowPloen1995,PloenPark1997,PloenPark1999,ParkKim2000}%
. It should be mentioned that the basic elements of a screw formulation for
MBS dynamics were already presented in \cite{Liu1988}, but did not receive
due attention.

A crucial feature of these geometric approaches is their frame invariance,
which allows for arbitrary representations of screws and for freely
assigning reference frames, which drastically simplify the kinematics
modeling and also provides a direct link to CAD models. Moreover, the POE,
and thus the kinematics, can even be formulated without the use of any joint
frame, which basically resembles the 'zero reference' formulation that was
reported for a robotic arm in \cite{Gupta1986}. On the other hand, classical
approaches to the description of joint kinematics are the Denavit-Hartenberg
(DH) \cite{DenavitHartenberg1955,KhalilKleinfinger1986} (in its different
forms) and the Sheth-Uicker two-frame convention \cite{ShethUicker1971}.
Such two-frame conventions are used in most of the current MBS dynamics
simulations packages that use relative coordinates. The Lie group
description, on the other hand, not only allows for arbitrary placement of
joint frames but makes them dispensable altogether.

The benefits of geometric modeling have been recognized already in robotics.
Recently, at least in robotics, the text books \cite%
{LynchPark2017,Murray,Selig} have reached a wider audience. Modern
approaches to robotics make extensive use of screw and Lie group theoretical
concepts. This is, also supported by the Universal Robot Description Format
(URDF) that is used for instance in the Robot Operating System (ROS), rather
than DH parameters. In MBS dynamics the benefits of geometric mechanics are
slowly being recognized. Interestingly, this mainly applies to the modeling
of MBS with flexible bodies undergoing large deformations \cite%
{Bauchau2011,Sonneville2014}. This is not surprising since geometrically
exact formulations require correct modeling of the finite kinematics of a
continua. The displacement field of a Cosserat beam, for instance, is a
proper motion in $E^{3}$ and thus modeled as motion in $SE\left( 3\right) $.
This is an extension of the original work on geometrically exact beams and
shells by Simo \cite{Simo1986,Simo1991} where the displacement field is
modeled on $SO\left( 3\right) \times {\mathbb{R}}^{3}$. Another topic where
Lie group theory is recently applied in MBS dynamics is the time
integration. To this end, Lie group integration schemes were modified and
applied to MBS models in absolute coordinate formulation \cite%
{BruelsCardonaArnold2012}, where the motion of individual bodies is
described as a general screw motion that are constrained according to the
interconnecting joints. It shall be remarked that, despite the current trend
to emphasize the use of Lie group (basic) concepts, the basics formulations
for non-linear flexible MBS were already reported by Borri et al., e.g. \cite%
{Borri2001a,Borri2001b,Borri2003}.

The aim of this paper is to provide a comprehensive summary of the basic
concepts for modeling MBS in terms of relative coordinates using joint
screws and to relate them to existing formulations that are scattered
throughout the literature. Without loss of generality the concepts are
introduced for a kinematic chain within a MBS with arbitrary topology \cite%
{Abhi2011a,Topology}. It is also the aim to show that MBS can be modeled in
a user-friendly way without having to follow restrictive modeling
conventions, and that this gives rise to $O\left( n\right) $ formulations.
The latter are not the topic of this paper.

The paper is organized as follows. In section \ref{secConfig} the MBS
configuration is described in terms of joint variables, used as generalized
coordinates, with the joint geometry parameterized by joint screw
coordinates. This classical approach of using body-fixed joint frames to
describe relative configurations is extended to a formulation that does not
involve joint frames. The corresponding relations for the MBS velocity are
derived in section \ref{secVel}. A formulation is introduced for each of the
four different definitions of rigid body twists found in the literature. The
latter are called the body-fixed, spatial, hybrid, and mixed twists. They
differ by the reference point used to measure the velocity and by the frame
in which the angular and translational velocities are resolved. The
different twist representations are introduced in appendix A.2. Recursive
relations for the respective Jacobians are derived, and the computational
aspects are discussed with emphasize on their decomposition. The presented
formulation allows for an efficient modeling of the MBS kinematics in terms
of readily available geometric data. Throughout the paper only a few basic
concepts from Lie group theory are required that are summarized in appendix
A. The used nomenclature is summarized in appendix B.

As for all Lie group formulations the biggest hurdle for a reader (who may
be already be familiar with MBS formulations) is the notation. The reader
not familiar with screws and Lie group modeling may want to consult the
appendix A.1 before reading section \ref{secConfig} and appendix A.2 before
reading section \ref{secVel}. This paper is aimed to provide a reference and
cannot replace an introductory textbook like \cite%
{LynchPark2017,Murray,Selig}. A beginner is recommended to consult \cite%
{LynchPark2017}. Yet there is no text book that treats the topic from a MBS
perspective. Readers not interested in the derivations could simply use the
main relations that are displayed with a black border.

\section{Configuration of a Kinematic Chain%
\label{secConfig}%
}

In this section the kinematics modeling using joint screw coordinates is
presented. For simplicity a single open kinematic chain is considered
comprising $n$ moving bodies interconnected by $n$ 1-DOF lower pair joints.
To simplify the formulation, but without loosing generality, higher-DOF
joints are modeled as combination of 1-DOF lower pair joints. Bodies and
joints are labeled with the same indices $i=1,\ldots ,n$ while the ground is
indexed with 0. With the sequential numbering of bodies and joints of the
kinematic chain, joint $i$ connects body $i$ to its predecessor body $i-1$.
A body-fixed reference frame (BFR) $\mathcal{F}_{i}$ is attached to body $i$
of the MBS. The body is then kinematically represented by this BFR. \vspace{%
-3ex}

\subsection{Joint Kinematics}

It has been the standard approach in MBS modeling to represent higher-DOF
joints by combination of 1-DOF lower pair joints, i.e. using either
revolute, prismatic, or screw joints. This will be adopted in the following
although this is not the way in which MBS models are implemented in
practice, but it simplifies the introduction of the presented approach
without compromising its generality. The justification of this approach is
that most technical joints are so-called lower kinematic pairs (also called
Reuleaux pairs) characterized by surface contact \cite%
{Reuleaux1875,Reuleaux1963}. That is, they are the mechanical generators of
motion subgroups of $SE\left( 3\right) $ \cite{Selig}. But not all motion
subgroups are generated by lower pairs. The 10 subgroups are listed in table %
\ref{TableSE3subgroups1}. So-called 'macro joints' are frequently used in
MBS modeling to generate motion subgroups by combination of lower pairs.
Table \ref{TableSE3subgroups2} shows the correspondence of motion subgroups
with lower pairs and macro joints. Missing in this list are joints relevant
for MBS modeling such as universal/hook and constant velocity joints since
they are not lower kinematic pairs. They can be modeled by combination of
lower pair joints.

\begin{table}[h] \centering%
\vspace{-3ex} 
\begin{tabular}[b]{cll}
\hline
$n$ & \textbf{Subgroup} & \textbf{Motion} \\ \hline
1 & $\mathbb{R}$ & 1-dim. translation along some axis%
\vspace{0.7ex}
\\ 
1 & $SO\left( 2\right) $ & 1-dim. rotation about arbitrary fixed axis%
\vspace{0.7ex}
\\ 
1 & $H_{p}$ & screw motion about arbitrary axis with finite pitch%
\vspace{0.7ex}
\\ 
2 & $\mathbb{R}^{2}$ & 2-dim. planar translation%
\vspace{0.7ex}
\\ 
2 & $SO\left( 2\right) \ltimes \mathbb{R}$ & translation along arbitrary
axis \& rotation along this axis%
\vspace{0.7ex}
\\ 
3 & $\mathbb{R}^{3}$ & spatial translations%
\vspace{0.7ex}
\\ 
3 & $SO\left( 3\right) $ & spatial rotations about arbitrary fixed point%
\vspace{0.7ex}
\\ 
3 & $H_{p}\ltimes \mathbb{R}^{2}$ & translation in a plane + screw motion $%
\bot $ to this plane (pitch $h$)%
\vspace{0.7ex}
\\ 
3 & $SO\left( 2\right) \ltimes \mathbb{R}^{2}=SE\left( 2\right) $ & planar
motions%
\vspace{0.7ex}
\\ 
4 & $SO\left( 2\right) \ltimes \mathbb{R}^{3}=SE\left( 2\right) \ltimes 
\mathbb{R}$ & spatial translations + rotation about axis with fixed
orientation \\ 
&  & (Sch\"{o}nflies motion)%
\vspace{0.7ex}
\\ 
6 & $SE\left( 3\right) $ & spatial motion \\ \hline
\end{tabular}%
\caption{$n$-dimensional motion subgroups of $SE(3)$.}\label%
{TableSE3subgroups1} \vspace{-3ex} 
\end{table}%
\begin{table}[h] \centering%
%
%
%
%
%
%
%
%
%
%
%
%
%
%
%
%
%
%
%
%
%
%
%
%
%
%
%
%
%
%
%
%
%
%
%
%
%
%
%
%
%
%
%
%
%
\begin{tabular}[b]{clll}
\hline
$n$ & \textbf{Subgroup} & \textbf{Lower Pair} & \textbf{Macro Joint} \\ 
\hline
1 & $\mathbb{R}$ & Prismatic Joint & \multicolumn{1}{c}{$\times $%
\vspace{0.7ex}%
} \\ 
1 & $SO\left( 2\right) $ & Revolute Joint & \multicolumn{1}{c}{$\times $%
\vspace{0.7ex}%
} \\ 
1 & $H_{p}$ & Screw Joint & \multicolumn{1}{c}{$\times $%
\vspace{0.7ex}%
} \\ 
2 & $\mathbb{R}^{2}$ & \multicolumn{1}{c}{$\times $%
\vspace{0.7ex}%
} & combination of two non-parallel \\ 
&  & \multicolumn{1}{c}{} & prismatic joints%
\vspace{0.7ex}
\\ 
2 & $SO\left( 2\right) \ltimes \mathbb{R}$ & Cylindrical Joint & 
\multicolumn{1}{c}{$\times $%
\vspace{0.7ex}%
} \\ 
3 & $\mathbb{R}^{3}$ & \multicolumn{1}{c}{$\times $%
\vspace{0.7ex}%
} & combination of three non-parallel \\ 
&  &  & prismatic joints%
\vspace{0.7ex}
\\ 
3 & $SO\left( 3\right) $ & Spherical Joint & \multicolumn{1}{c}{$\times $%
\vspace{0.7ex}%
} \\ 
3 & $H_{p}\ltimes \mathbb{R}^{2}$ & \multicolumn{1}{c}{$\times $%
\vspace{0.7ex}%
} & planar joint + screw joint with axis \\ 
&  &  & normal to plane \\ 
3 & $SO\left( 2\right) \ltimes \mathbb{R}^{2}=SE\left( 2\right) $ & Planar
Joint & \multicolumn{1}{c}{$\times $%
\vspace{0.7ex}%
} \\ 
4 & $SO\left( 2\right) \ltimes \mathbb{R}^{3}=SE\left( 2\right) \ltimes 
\mathbb{R}$ & \multicolumn{1}{c}{$\times $%
\vspace{0.7ex}%
} & planar joint + prismatic joint \\ 
&  & \multicolumn{1}{c}{} & with axis normal to plane%
\vspace{0.7ex}
\\ 
6 & $SE\left( 3\right) $ & \multicolumn{1}{c}{$\times $%
\vspace{0.7ex}%
} & 'free joint' \\ \hline
\end{tabular}%
\caption{Mechanical generators of the $n$-dimensional subgroups of $SE(3)$. 
A motion subgroup can be generated by a lower pair or by a 'macro joint', i.e. a combination of joints with smaller DOF.}%
\label{TableSE3subgroups2}%
\end{table}%

\begin{figure}[t]
\centerline{
\includegraphics[width=11cm]{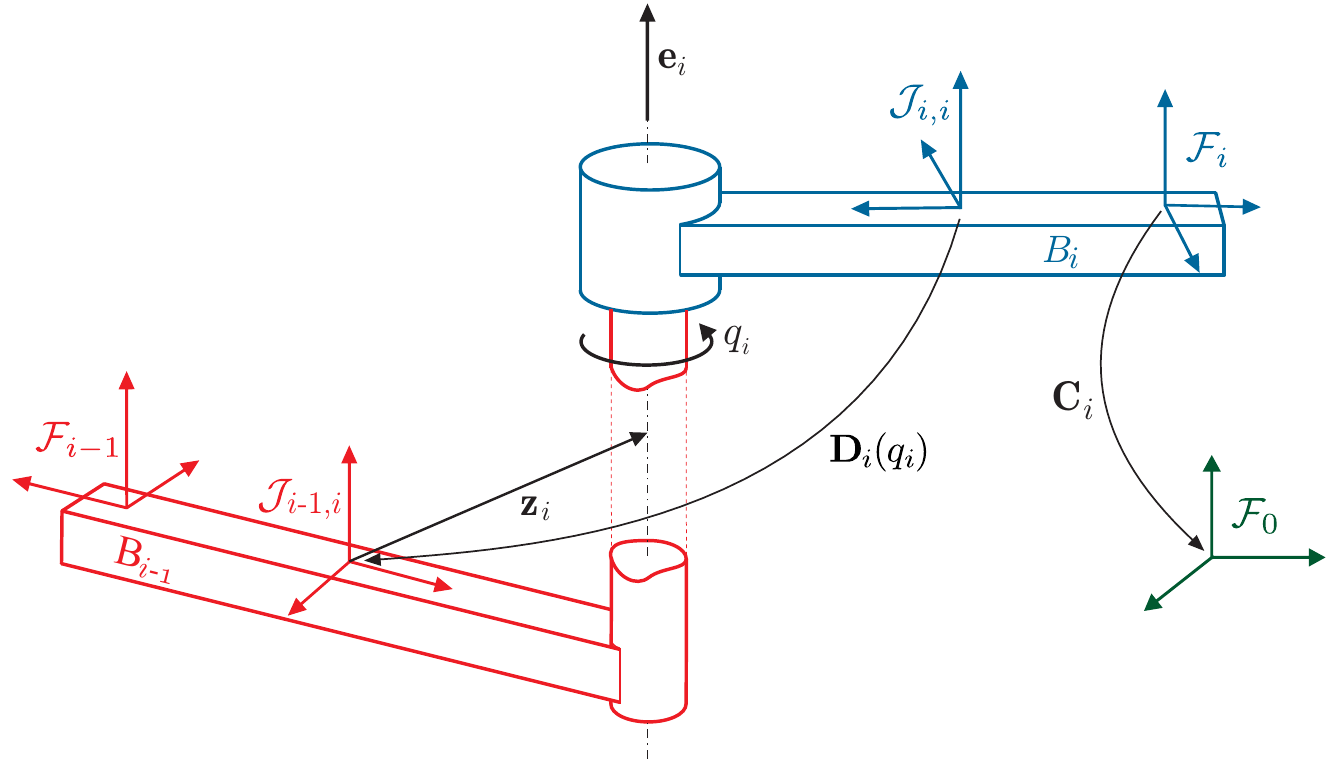}}
\caption{Description of the kinematics of joint $i$ connecting body $i$ with
its predecessor body $i-1$. A body-fixed JFR $\mathcal{J}_{i,i}$ is
introduced on body $i$, and $\mathcal{J}_{i-1,i}$, on body $i-1$,
respectively. A revolute joint is shown as example.}
\label{figJointKinematics_Z}
\end{figure}

The classical approach to describe joint kinematics is to introduce an
additional pair of body-fixed \emph{joint frames} (JFR) for each joint (fig. %
\ref{figJointKinematics_Z}) \cite{Uicker2013}. Denote with $\mathcal{J}%
_{i-1,i}$ the JFR for joint $i$ on body $i-1$ and with $\mathcal{J}_{i,i}$
the JFR on body $i$. The relative motion of adjacent bodies is represented
by the frame transformation between the respective JFRs that can be
described in terms of screw coordinates (appendix A.1). 

Lower pair 1-DOF joints restrict the interconnected bodies so to perform
screw motions with a certain pitch $h$. Revolute joints have pitch $h=0$,
and prismatic joints $h=\infty $, while proper screw joints have a non-zero
finite pitch. Denote with%
\begin{equation}
{^{i-1}}\mathbf{Z}_{i}=\left( 
\begin{array}{c}
{^{i-1}}\mathbf{e}{_{i}} \\ 
{^{i-1}}\mathbf{z}{_{i}}\times {^{i-1}}\mathbf{e}{_{i}}+{^{i-1}}\mathbf{e}{%
_{i}h}_{i}%
\end{array}%
\right)  \label{Zi}
\end{equation}%
the unit screw coordinate vector of joint $i$ expressed in the JFR $\mathcal{%
J}_{i-1,i}$ on body $i-1$, where ${^{i-1}}\mathbf{z}{_{i}}$ is the position
vector of a point on the joint axis measured in the JFR $\mathcal{J}_{i-1,i}$%
, and ${^{i-1}}\mathbf{e}{_{i}}$ is the unit vector along the joint axis
resolved in JFR $\mathcal{J}_{i-1,i}$.

\begin{assumption}
It is assumed throughout the paper that the two JFRs coincide in the
reference configuration $q_{i}=0$. This assumption can be easily relaxed if
required.
\end{assumption}

Denote with $q_{i}$ the joint variable (angle, translation). With the above
assumption the configuration of the JFR $\mathcal{J}_{i,i}$ on body $i$
relative to the JFR $\mathcal{J}_{i-1,i}$ on body $i-1$ is given by the
exponential in (\ref{SE3expP}) as $\mathbf{D}_{i}\left( q_{i}\right) :=\exp
({}{}{^{i-1}}\mathbf{Z}_{i}q_{i})$.

\begin{remark}
\label{assumeJFR}%
It is common practice to locate the JFRs with their origins at the joint
axis (as in Fig. 2), so that $\mathbf{z}=\mathbf{0}$. Then the joint screw
coordinates for the three types of 1-DOF joints are%
\begin{equation}
\mathbf{Z}^{\text{revolute}}=\left( 
\begin{array}{c}
\mathbf{e} \\ 
\mathbf{0}%
\end{array}%
\right) ,\ \mathbf{Z}^{\text{screw}}=\left( 
\begin{array}{c}
\mathbf{e} \\ 
\mathbf{e}{h}%
\end{array}%
\right) ,\ \mathbf{Z}^{\text{prismatic}}=\left( 
\begin{array}{c}
\mathbf{0} \\ 
\mathbf{e}%
\end{array}%
\right) .
\end{equation}
\end{remark}

\subsection{Recursive Kinematics using Body-Fixed Joint Frames}

The absolute configuration of body $i$, i.e. the configuration of its BFR $%
\mathcal{F}_{i}$ relative to the inertial frame (IFR) $\mathcal{F}_{0}$, is
represented by $\mathbf{C}_{i}\in SE\left( 3\right) $. The \emph{relative
configuration} of body $i$ relative to body $i-1$ is $\mathbf{C}_{i-1,i}:=%
\mathbf{C}_{i-1}^{-1}\mathbf{C}_{i}$. The configuration of a rigid body in
the kinematic chain can be determined recursively by successive combination
of the relative configurations of adjacent bodies as $\mathbf{C}_{i}=\mathbf{%
C}_{0,1}\mathbf{C}_{1,2}\cdots \mathbf{C}_{i-1,i}$.

For joint $i$ denote with ${}\mathbf{S}_{i-1,i}$ the constant transformation
from JFR $\mathcal{J}_{i-1,i}$ to the RFR $\mathcal{F}_{i-1}$ on body $i-1$,
and with $\mathbf{S}_{i,i}$ the constant transformation from JFR $\mathcal{J}%
_{i,i}$ to the RFR $\mathcal{F}_{i}$ on body $i$ (fig. \ref%
{figJointKinematics_JFR}). The relative configuration is then $\mathbf{C}%
_{i-1,i}={}\mathbf{S}_{i-1,i}\mathbf{D}_{i}\left( q_{i}\right) {}\mathbf{S}%
_{i,i}^{-1}$. Denote with $\mathbf{q}\in {\mathbb{V}}^{n}$ the vector of
joint variables that serve as generalized coordinates of the MBS. The joint
space manifold is ${\mathbb{V}}^{n}={\mathbb{R}}^{n_{\text{P}}}\times 
\mathbb{T}^{n_{\text{R}}}$ for an MBS model comprising $n_{\text{P}}$
prismatic and $n_{\text{R}}$ revolute/screw joints ($n_{\text{P}}+n_{\text{R}%
}=n$).

The \emph{absolute configuration} (i.e. relative to the IFR) of body $i$ in
the chain is%
\begin{equation}
\begin{tabular}{|lll|}
\hline
&  &  \\ 
& $\mathbf{C}_{i}\left( \mathbf{q}\right) =\mathbf{S}_{0,1}\mathbf{D}%
_{1}\left( q_{1}\right) \mathbf{S}_{1,1}^{-1}\cdot \mathbf{S}_{1,2}\mathbf{D}%
_{2}\left( q_{2}\right) \mathbf{S}_{2,2}^{-1}\cdot \ldots \cdot \mathbf{S}%
_{i-1,i}\mathbf{D}_{i}\left( q_{i}\right) \mathbf{S}_{i,i}^{-1}.$ &  \\ 
&  &  \\ \hline
\end{tabular}
\label{recConf}
\end{equation}

This formulation requires the following modeling steps:

\begin{itemize}
\item Introduction of body-fixed JFR $\mathcal{J}_{i,i}$ at body $i$ with
relative configuration $\mathbf{S}_{i,i}$,

\item Introduction of body-fixed JFR $\mathcal{J}_{i-1,i}$ at body $i-1$
with relative configurations $\mathbf{S}_{i-1,i}$,

\item The screw coordinate vector ${^{i-1}}\mathbf{Z}_{i}$ of joint $i$
represented in JFR $\mathcal{J}_{i-1,i}$ at body $i-1$.
\end{itemize}

The expression (\ref{recConf}) is the standard MBS formulation for the
kinematics of an open chain in terms of \emph{relative coordinates}, i.e.
joint angles or translations. For 1-DOF joints the JFR is usually oriented
such that its 3-axis points along the joint axis (as in figure \ref%
{figJointKinematics_JFR}). The screw coordinates are then ${^{i-1}}\mathbf{Z}%
_{i}=(0,0,1-s_{i},0,0,s_{i}+h_{i}\left( 1-s_{i}\right) )^{T}$ where $%
s_{i}\,=1$ for prismatic joint, and $s_{i}=0$ for a screw joint with finite
pitch $h_{i}$ (for revolute joints $h_{i}=0$). 
\begin{figure}[b]
\centerline{
\includegraphics[width=14cm]{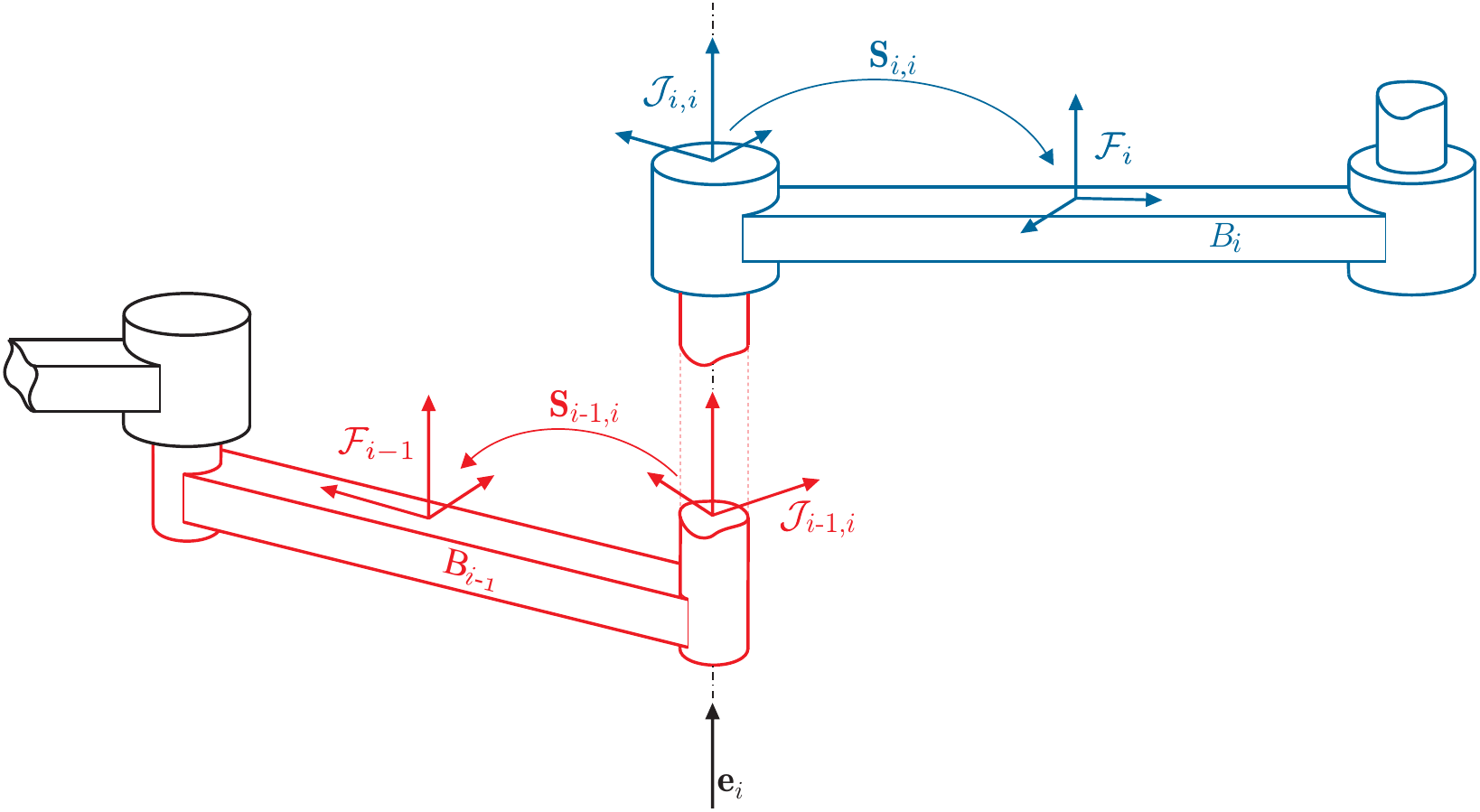}}
\caption{Definition of body-fixed RFR $\mathcal{F}_{i}$ and JFR $\mathcal{J}%
_{i,i}$ and $\mathcal{J}_{i-1,i}$, for joint $i$.}
\label{figJointKinematics_JFR}
\end{figure}

\begin{remark}
The matrix $\mathbf{C}_{i}$ is used to represent the configuration of body $%
i $; hence the symbol. Frequently the symbol $\mathbf{T}_{i}$ is used \cite%
{LynchPark2017,Uicker2013}, which refers to the fact that these matrices
describe the transformation of point coordinates (appendix A.1)
\end{remark}

\begin{remark}
It is important to emphasize that the Lie group formulation (\ref{recConf})
is merely another approach to the standard matrix formulation of MBS
kinematics aiming at compact expressions that simplify the implementation
without compromising the efficiency. It also includes the various
conventions used to describe the joint kinematics. An excellent overview of
classical matrix methods (also with emphasis on how they can be employed for
synthesis) can be found in \cite{Uicker2013}. For instance, $\mathbf{S}%
_{i-1,i}$ and $\mathbf{S}_{i,i}$ can be parameterized in terms of the
constant part of the Denavit-Hartenberg (DH) parameters \cite{Uicker2013}.
The formulation (\ref{recConf}) in particular resembles the Sheth-Uicker
convention (that was introduced to eliminate the ambiguity of DH parameter) 
\cite{ShethUicker1971,Uicker2013}. In that notation the matrices $\mathbf{S}%
_{i-1,i}$ and $\mathbf{S}_{i,i}$ are called the \emph{shape matrices} of
joint $i$. However, the Sheth-Uicker convention still presumes certain
alignment of joint axes. E.g. a revolute axis is supposed to be parallel to
the 3-axis of the JFRs. A recent discussion of these notations can be found
in \cite{BongardtPhD}. An expression similar to (\ref{recConf}) was also
presented in \cite{Orin1979} where no restriction on the joint axis is
imposed. A recursive formulation of the MBS motion equations using
homogeneous transformation matrices was also presented in \cite%
{Legnani1,Legnani2}.
\end{remark}

\begin{remark}[Multi-DOF joints]
The description for 1-DOF joints in terms of a screw coordinate vector $%
\mathbf{Z}_{i}$ can be generalized to joints with more than one DOF. For a
joint with DOF $\nu $ the relative configuration of the JFRs can
alternatively be described in terms of $\nu $ joint variables $%
q_{i_{1}},\ldots ,q_{i_{\nu }}$ by $\mathbf{D}_{i}\left( q_{i_{1}},\ldots
,q_{i_{\nu }}\right) :=\exp ({}{}{^{i-1}}\mathbf{Z}_{i_{1}}q_{i_{1}}+\ldots +%
{^{i-1}}\mathbf{Z}_{i_{\nu }}q_{i_{\nu }})$ or $\mathbf{D}_{i}\left(
q_{i_{1}},\ldots ,q_{i_{\nu }}\right) :=\exp ({}{}{^{i-1}}\mathbf{Z}%
_{i_{1}}q_{i_{1}})\cdot \ldots \cdot \exp ({}{}{^{i-1}}\mathbf{Z}_{i_{\nu
}}q_{i_{\nu }})$. For a spherical joint, for instance, the variables in the
first form are the components of the rotation axis times angle in (\ref%
{SO3exp}), and in the second form these are three angles corresponding to
the order of 1-DOF rotations (e.g. Euler-angles). For lower pair joints, in
the first case, $q_{i_{1}},\ldots ,q_{i_{\nu }}$ are canonical coordinates
of first kind on the joint's motion subgroup, and in the second case they
are canonical coordinates of second kind \cite{Murray}. The $\mathbf{Z}%
_{i_{1}},\ldots ,\mathbf{Z}_{i_{\nu }}$ form a basis on the subalgebra of
the motion subgroup generated by the joint.
\end{remark}

\subsection{Recursive Kinematics without Body-Fixed Joint Frames}

The introduction of joint frames is a tedious step within the MBS kinematics
modeling. Moreover, it is desirable to minimize the data required to
formulate the kinematic relations. In this regard the frame invariance of
screws is beneficial.

The two constant transformations from the JFR to the BFR on the respective
body can be summarized using (\ref{Adexp}) as%
\begin{equation}
\mathbf{C}_{i-1,i}\left( q_{i}\right) =\mathbf{S}_{i-1,i}\mathbf{D}%
_{i}\left( q_{i}\right) \mathbf{S}_{i,i}^{-1}=\mathbf{S}_{i-1,1}\mathbf{S}%
_{i,i}^{-1}\mathbf{S}_{i,i}\mathbf{D}_{i}\left( q_{i}\right) \mathbf{S}%
_{i,i}^{-1}=\mathbf{B}_{i}\exp ({^{i}}\mathbf{X}_{i}q_{i})  \label{relConf}
\end{equation}%
so that the relative configuration splits into only one constant and a
variable part. The constant part $\mathbf{B}_{i}:=\mathbf{S}_{i-1,i}\mathbf{S%
}_{i,i}^{-1}=\mathbf{C}_{i-1,i}\left( 0\right) $ is the reference
configuration of body $i$ w.r.t. body $i-1$ when $q_{i}=0$. The variable
part is now given in terms of the constant screw coordinate vector of joint $%
i$ represented in BFR $\mathcal{F}_{i}$ 
\begin{equation}
{^{i}}\mathbf{X}_{i}=\mathbf{Ad}_{\mathbf{S}_{i,i}}{^{i-1}}\mathbf{Z}%
_{i}=\left( 
\begin{array}{c}
{^{i}}\mathbf{e}_{i} \\ 
{^{i}}\mathbf{x}_{i,i}\times {^{i}}\mathbf{e}_{i}+h_{i}{^{i}}\mathbf{e}_{i}%
\end{array}%
\right) .  \label{Xb}
\end{equation}%
The matrix $\mathbf{Ad}_{\mathbf{S}_{i,i}}$, defined in (\ref{Ad}),
transforms screw coordinates represented in $\mathcal{J}_{i,i-1}$ to those
represented in $\mathcal{F}_{i}$ according to their relative configuration
described by $\mathbf{S}_{i,i}$.

\begin{figure}[b]
\centerline{
\includegraphics[width=14cm]{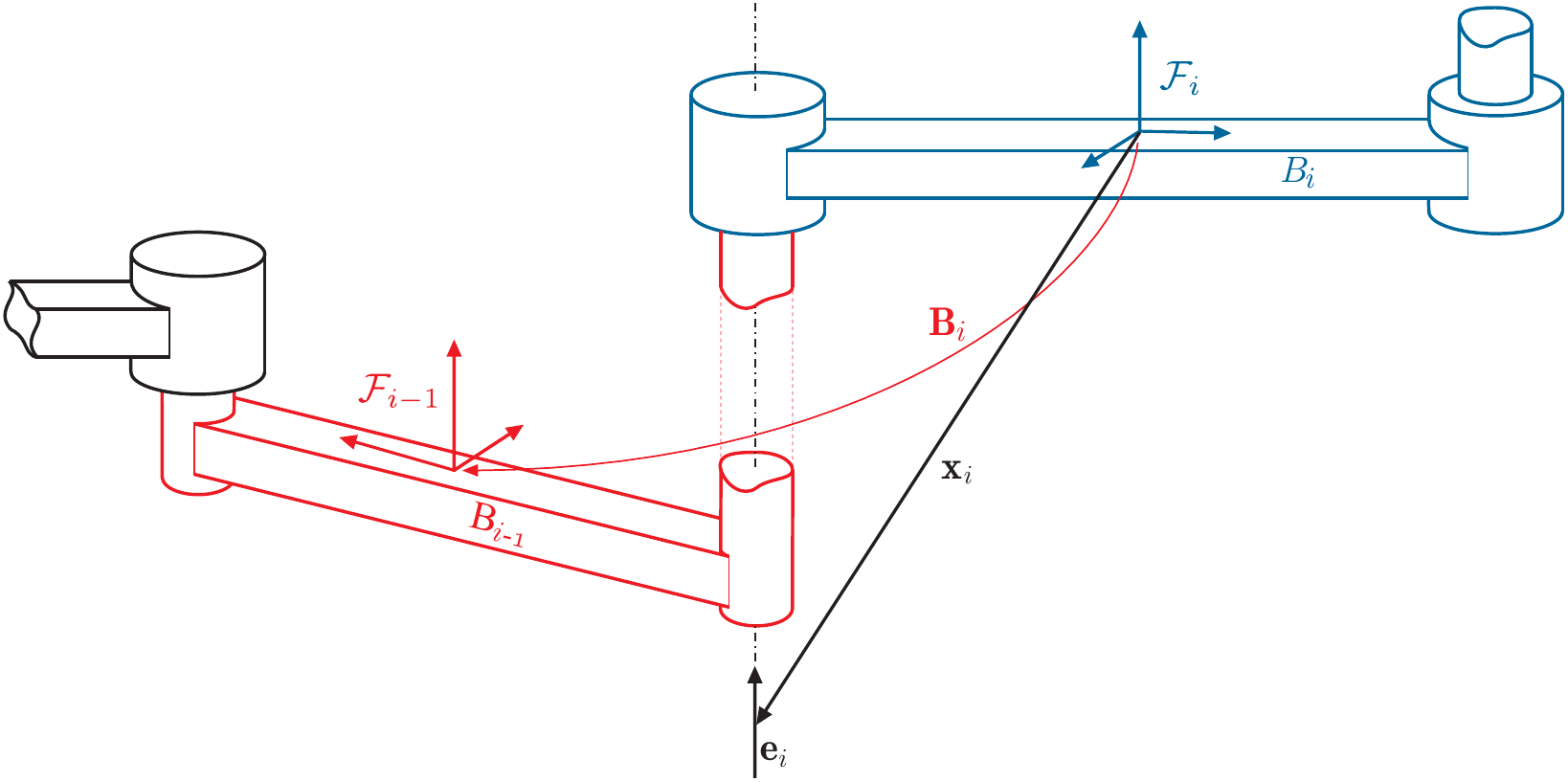}}
\caption{Description of the kinematics of joint $i$ without body-fixed JFRs,
in the zero-reference relative configuration with $q_{i}=0$. The vector $%
\mathbf{x}_{i}$ is used when the joint screw coordinates are represented in
the BFR $\mathcal{F}_{i}$ on body $i$, and $\bar{\mathbf{x}}_{i-1,i}$ is
used when the joint screw coordinates are represented in the BFR $\mathcal{F}%
_{i-1}$ on body $i-1$.}
\label{figJointKinematics_BFR}
\end{figure}

As indicated in fig. \ref{figJointKinematics_BFR}, here ${^{i}}\mathbf{e}%
_{i} $ is the unit vector along the axis of joint $i$ resolved in the BFR $%
\mathcal{F}_{i}$, and ${^{i}}\mathbf{x}_{i,i}$ is the position vector of a
point on the axis of joint $i$, measured and resolved in $\mathcal{F}_{i}$.
This is indeed the same screw as in (\ref{Zi}) but expressed in the BFR on
body $i$. 

The joint screw can alternatively be represented in $\mathcal{F}_{i-1}$.
Then 
\begin{equation}
\mathbf{C}_{i-1,i}\left( q_{i}\right) =\mathbf{S}_{i,i-1}\mathbf{D}%
_{i}\left( q_{i}\right) \mathbf{S}_{i,i}^{-1}=\mathbf{S}_{i,i-1}\mathbf{D}%
_{i}\left( q_{i}\right) \mathbf{S}_{i,i-1}^{-1}\mathbf{S}_{i,i-1}\mathbf{S}%
_{i,i}^{-1}=\exp ({^{i-1}}\bar{\mathbf{X}}_{i}q_{i})\mathbf{B}_{i}
\end{equation}%
with the joint screw coordinate vector%
\begin{equation}
{^{i-1}}\bar{\mathbf{X}}_{i}=\mathbf{Ad}_{\mathbf{S}_{i-1,i}}{^{i-1}}\mathbf{%
Z}_{i}=\mathbf{Ad}_{\mathbf{B}_{i}}{^{i}}\mathbf{X}_{i}=\left( 
\begin{array}{c}
{^{i-1}}\mathbf{e}_{i} \\ 
{^{i-1}}\bar{\mathbf{x}}_{i-1,i}\times {^{i-1}}\mathbf{e}_{i}+{^{i-1}}%
\mathbf{e}_{i}h_{i}%
\end{array}%
\right)  \label{Xibar}
\end{equation}%
now expressed in the BFR $\mathcal{F}_{i-1}$ at body $i-1$, where ${^{i-1}}%
\bar{\mathbf{x}}_{i-1,i}$ is the position vector of a point on the axis of
joint $i$ measured in $\mathcal{F}_{i-1}$.

Successive combination of the relative configurations yields%
\begin{equation}
\begin{tabular}{|llll|}
\hline
&  &  &  \\ 
& $\mathbf{C}_{i}\left( \mathbf{q}\right) 
\hspace{-2ex}%
$ & $=\mathbf{B}_{1}\exp ({^{1}}\mathbf{X}_{1}q_{1})\cdot \mathbf{B}_{2}\exp
({^{2}}\mathbf{X}_{2}q_{2})\cdot \ldots \cdot \mathbf{B}_{i}\exp ({^{i}}%
\mathbf{X}_{i}q_{i})$ &  \\ 
&  &  & $%
\vspace{-2ex}%
$ \\ 
&  & $=\exp ({^{0}}\bar{\mathbf{X}}_{1}q_{1})\mathbf{B}_{1}\cdot \exp ({^{1}}%
\bar{\mathbf{X}}_{2}{q}_{2})\mathbf{B}_{2}\cdot \ldots \cdot \exp ({^{i-1}}%
\bar{\mathbf{X}}_{i}q_{i})\mathbf{B}_{i}.$ &  \\ 
&  &  &  \\ \hline
\end{tabular}
\label{POEX}
\end{equation}%
The first form of (\ref{POEX}) was reported \cite{PloenPark1999}, and both
forms in \cite{Park1994,ParkBobrowPloen1995}. It will be called the \emph{%
body-fixed Product-of-Exponentials} (POE) formula in body-fixed description
since the joint kinematics is expressed by exponentials of joint screws. It
seems to be more convenient to work with the screw coordinates ${^{i}}%
\mathbf{X}_{i}$. Also in \cite{Angeles2003} two variants of the kinematic
description of a serial chain were presented using a BFR on body $i-1$ or $i$%
, respectively.

In summary this body-fixed POE formulation does not require introduction of
JFRs. It only requires the following readily available information:

\begin{itemize}
\item The relative reference configuration $\mathbf{B}_{i}$ of the adjacent
bodies connected by joint $i$ for $q_{i}=0$,

\item The screw coordinates ${^{i}}\mathbf{X}_{i}$ of joint $i$ represented
in the BFR $\mathcal{F}_{i}$ at body $i$, or alternatively the screw
coordinates ${^{i-1}}\bar{\mathbf{X}}_{i}$ represented in the BFR $\mathcal{F%
}_{i-1}$ at body $i-1$.
\end{itemize}

The form (\ref{POEX}) simplifies the expression for the joint kinematics.
Its main advantage is that it only involves the reference configuration $%
\mathbf{B}_{i}$ of BFRs.

\subsection{Recursive Kinematics without Body-Fixed Joint Frames and Screw
Coordinates}

\begin{figure}[bh]
\centerline{
\includegraphics[width=13cm]{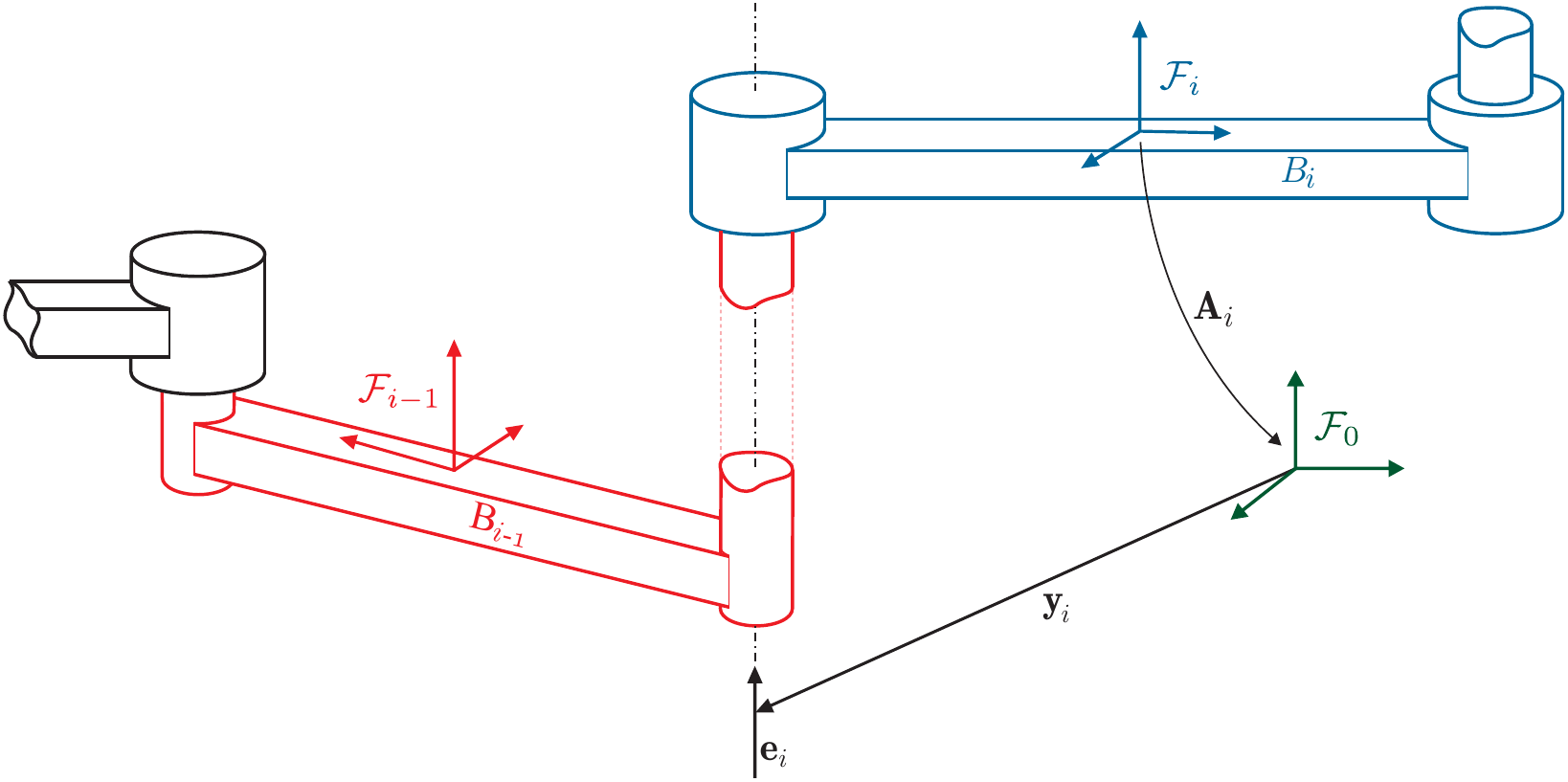}}
\caption{Description of the kinematics of joint $i$ with respect to the
spatial IFR in the zero-reference configuration with $\mathbf{q}=\mathbf{0}$%
. }
\label{figJointKinematics_IFR}
\end{figure}

Thanks to the frame invariance, the joint screw coordinates can even be
described in the spatial IFR, i.e. without reference to any body-fixed
frames. To this end, (\ref{POEX}) is written as%
\begin{eqnarray}
\mathbf{C}_{i}\left( \mathbf{q}\right) &=&\mathbf{B}_{1}\exp ({^{1}}\mathbf{X%
}_{1}q_{1})\mathbf{B}_{1}^{-1}\cdot  \notag \\
&&\mathbf{B}_{1}\mathbf{B}_{2}\exp ({^{2}}\mathbf{X}_{2}q_{2})\mathbf{B}%
_{2}^{-1}\mathbf{B}_{1}^{-1}\cdot \ldots \\
&&\mathbf{B}_{1}\cdots \mathbf{B}_{i}\exp ({^{i}}\mathbf{X}_{i}q_{i})\mathbf{%
B}_{i}^{-1}\cdots \mathbf{B}_{1}^{-1}\mathbf{B}_{1}\cdots \mathbf{B}_{i}. 
\notag
\end{eqnarray}%
Relation (\ref{Adexp}) yields $\mathbf{B}\exp (q\widehat{\mathbf{X}})\mathbf{%
B}^{-1}=\exp (q\mathbf{B}\widehat{\mathbf{X}}\mathbf{B}^{-1})=$ $\exp (q%
\mathbf{Ad}_{\mathbf{B}}\mathbf{X})$ so that%
\begin{equation}
\begin{tabular}{|lll|}
\hline
&  &  \\ 
& $\mathbf{C}_{i}\left( \mathbf{q}\right) =\exp (\mathbf{Y}_{1}q_{1})\cdot
\exp (\mathbf{Y}_{2}q_{2})\cdot \ldots \cdot \exp (\mathbf{Y}_{i}q_{i})%
\mathbf{A}_{i}.$ &  \\ 
&  &  \\ \hline
\end{tabular}
\label{POEY}
\end{equation}%
Here 
\begin{equation}
\mathbf{A}_{i}:=\mathbf{B}_{1}\cdots \mathbf{B}_{i}=\mathbf{C}_{i}\left( 
\mathbf{0}\right) =\left( 
\begin{array}{cc}
\mathbf{R}_{i}\left( \mathbf{0}\right) & \mathbf{r}_{i}\left( \mathbf{0}%
\right) \\ 
\mathbf{0} & 1%
\end{array}%
\right)  \label{Ai}
\end{equation}%
is the absolute reference configuration (i.e. relative to IFR) of body $i$,
and%
\begin{equation}
\mathbf{Y}_{j}=\mathbf{Ad}_{\mathbf{A}_{j}}{^{j}}\mathbf{X}_{j}=\left( 
\begin{array}{c}
\mathbf{e}_{j} \\ 
\mathbf{y}_{j}\times \mathbf{e}_{j}+h_{j}\mathbf{e}_{j}%
\end{array}%
\right)  \label{XY}
\end{equation}%
is the screw coordinate vector of joint $j$ represented in the IFR $\mathcal{%
F}_{0}$ in the reference configuration with $\mathbf{q}=\mathbf{0}$ (fig. %
\ref{figJointKinematics_IFR}). The direction unit vector $\mathbf{e}_{j}$
and the position vector $\mathbf{y}_{j}$ of a point on the joint axis are
expressed in the IFR $\mathcal{F}_{0}$ (leading superscript '0' omitted).
The transformation (\ref{XY}) relates the body-fixed to the spatial
representation of joint screw in the reference configuration. The product of
the exp mappings in (\ref{POEY}) describes the motion of a RFR on body $i$,
which at $\mathbf{q}=\mathbf{0}$ coincides with the IFR, relative to the
IFR. The relation to the actual BFR is achieved by the subsequent
transformation $\mathbf{A}_{i}$. Such a 'zero reference' formulation has
been first reported by Gupta \cite{Gupta1986} in terms of frame
transformation matrices, and was latter introduced by Brockett \cite%
{Brockett1984} as the POE formula for robotic manipulators. The formulation (%
\ref{POEY}) was then used in \cite{BrockettStokesPark1993} for MBS modeling.
It should be remarked that in the classical literature on screws, the
spatial representation of a screw is denoted with the '\$' symbol \cite%
{HuntBook1990,Hunt1991}.

All data required for this \emph{spatial POE formulation} is represented in
the spatial IFR:

\begin{itemize}
\item Absolute reference configurations $\mathbf{A}_{i}=\mathbf{C}_{i}\left( 
\mathbf{0}\right) $, i.e. the reference configuration of body $i$ w.r.t. the
IFR $\mathcal{F}_{0}$ for $\mathbf{q}=\mathbf{0}$,

\item Joint screw coordinates $\mathbf{Y}_{i}\equiv $ ${^{0}}\mathbf{Y}%
_{i}^{0}$ in spatial representations, i.e. measured and resolved in the IFR $%
\mathcal{F}_{0}$ for $\mathbf{q}=\mathbf{0}$.
\end{itemize}

The result (\ref{POEY}) is remarkable since it allows for formulating the
MBS kinematics \emph{without body-fixed joint frames}. From a modeling
perspective this has proven very useful since no joint transformations $%
\mathbf{S}_{i,i},\mathbf{S}_{i-1,i}$ or $\mathbf{B}_{i}$ are needed. Only
required are the absolute reference configurations $\mathbf{A}_{i}$ w.r.t.
to the IFR, and the reference screw coordinates (\ref{XY}), i.e. $\mathbf{e}%
_{i}$ and $\mathbf{p}_{i}$, resolved in the IFR. This is in particular
advantageous when processing CAD data. Moreover, if in the reference
(construction) configuration the RFR of the bodies coincide with the IFR
(global CAD reference system), i.e. all parts are designed w.r.t. the same
RFR, then $\mathbf{A}_{i}=\mathbf{I}$ and $\mathbf{Y}_{j}={^{j}}\mathbf{X}%
_{j}$. 

\subsection{Example%
\label{secExamp1}%
}

\begin{figure}[b]
\centerline{\includegraphics[width=10cm]{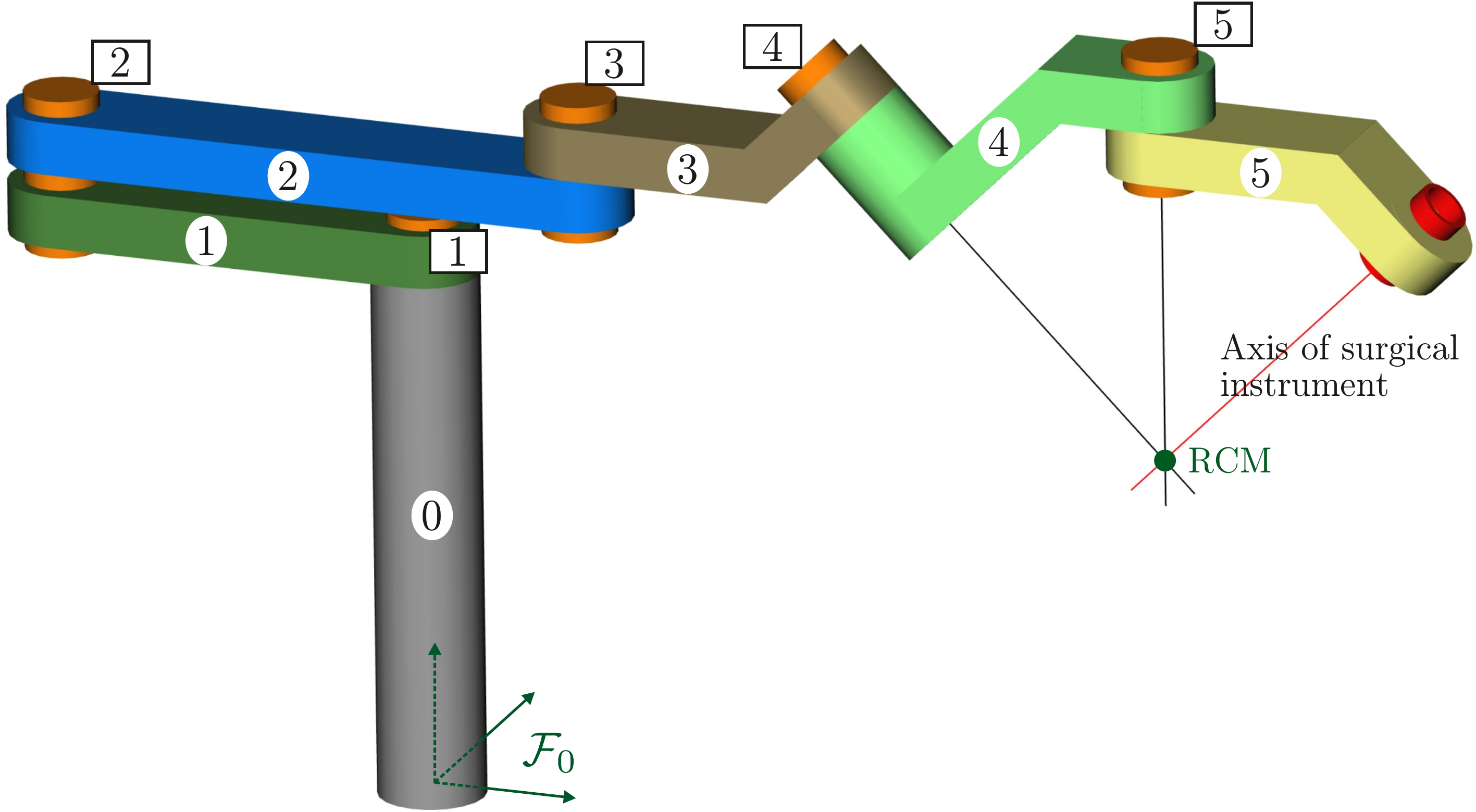}}
\caption{Model of the RCM mechanism disclosed in \protect\cite%
{WonChoiPeine2009}. The model was create with the MBS tool Alaska.}
\label{figRCMRef}
\end{figure}

\begin{figure}[b]
\centerline{\includegraphics[width=13cm]{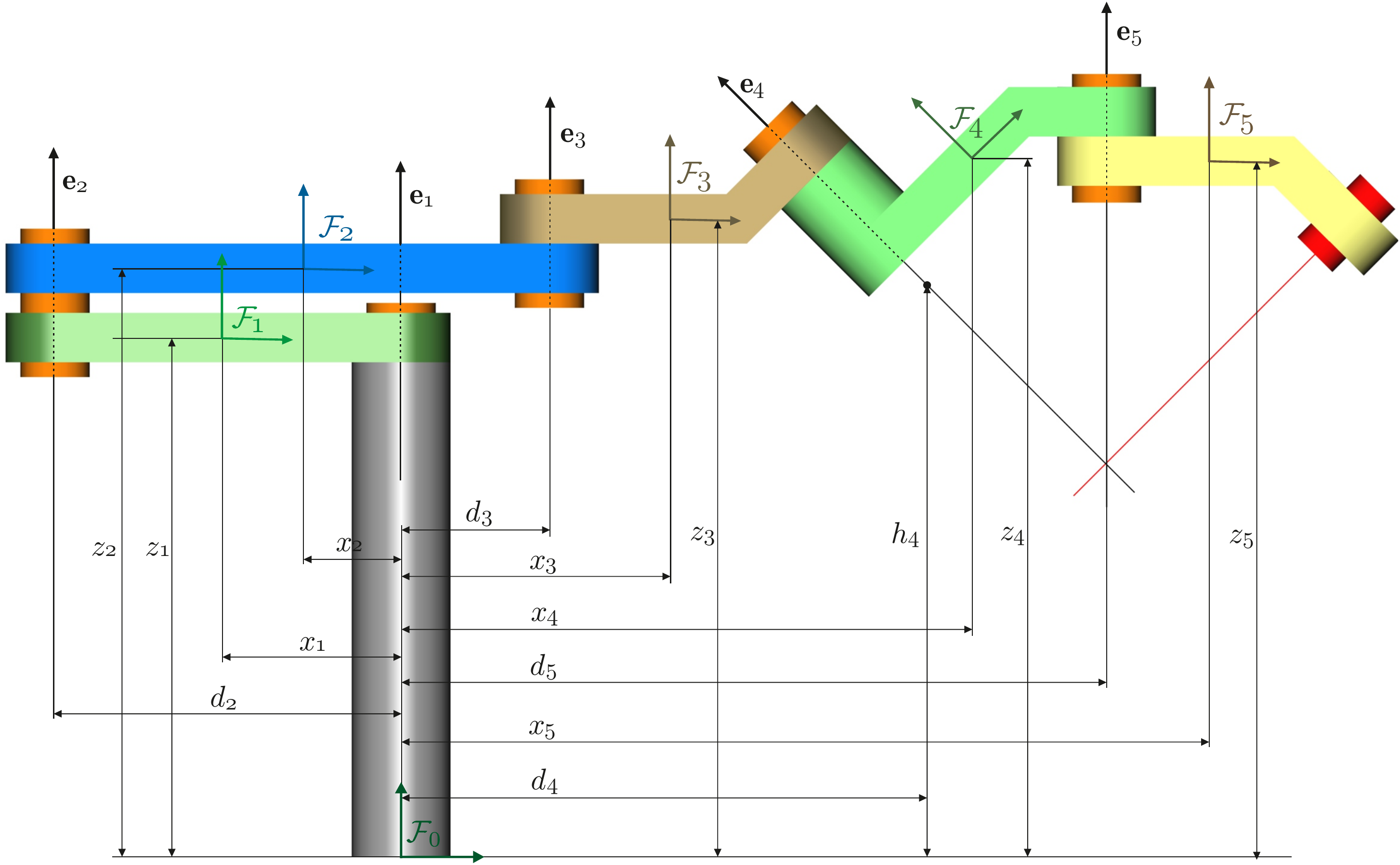}}
\caption{Description of the geometry of the RCM mechanism.}
\label{figRCMGeom}
\end{figure}

Figure \ref{figRCMRef} shows a surgical device that consists of a robot arm
and a remote center of motion (RCM) mechanism. This was disclosed in the
patent \cite{WonChoiPeine2009}. The robot arm consisting of bodies 1,2,3 is
used to position the RCM mechanism consisting of the bodies 4 and 5. The
surgical instrument is mounted in the socket at the remote end of body 5.
The axes of joints 4 and 5 and of the instrument intersect at one point.
This allows the instrument to freely pivot around an incision point.

The reference configuration is shown in fig. \ref{figRCMRef}. The IFR is
located at the base of the mechanism. The joint screw coordinates in spatial
representation are determined by the geometric parameters shown in figure %
\ref{figRCMGeom}. The position vectors $\mathbf{y}_{i}$ and unit vectors $%
\mathbf{e}_{i}$ in (\ref{XY}) are 
\begin{eqnarray*}
\mathbf{y}_{1} &=&\left( 0,0,0\right) ^{T},\mathbf{y}_{2}=\left(
-d_{2},0,0\right) ^{T},\mathbf{y}_{3}=\left( d_{3},0,0\right) ^{T},\mathbf{y}%
_{4}=\left( d_{4},0,h_{4}\right) ^{T},\mathbf{y}_{5}=\left( d_{5},0,0\right)
^{T} \\
\mathbf{e}_{1} &=&\mathbf{e}_{2}=\mathbf{e}_{3}=\mathbf{e}_{5}=\left(
0,0,1\right) ^{T},\mathbf{e}_{4}=(-(1/\sqrt{2},0,1/\sqrt{2})^{T}.
\end{eqnarray*}%
Since any point on the joint axes can be used, the 3-components in $\mathbf{y%
}_{i},i=1,2,3,5$ are set to zero. An arbitrary point on the axis of joint 4
is chosen as indicated. The joint screw coordinates (\ref{XY}) are thus%
\begin{eqnarray*}
\mathbf{Y}_{1} &=&\left( 0,0,1,0,0,0\right) ^{T},\mathbf{Y}_{2}=\left(
\{0,0,1,0,d_{2},0\right) ^{T},\mathbf{Y}_{3}=\left( 0,0,1,0,-d_{3},0\right)
^{T} \\
\mathbf{Y}_{4} &=&(-1\sqrt{2},0,1/\sqrt{2},0,-d_{4}/\sqrt{2}-h_{4}/\sqrt{2}%
,0)^{T},\mathbf{Y}_{5}=\left( 0,0,1,0,-d_{5},0\right) ^{T}.
\end{eqnarray*}%
The reference configurations (\ref{Ai}) of the bodies are determined by%
\begin{eqnarray*}
\mathbf{R}_{1}\left( \mathbf{0}\right) &=&\mathbf{R}_{2}\left( \mathbf{0}%
\right) =\mathbf{R}_{3}\left( \mathbf{0}\right) =\mathbf{R}_{5}\left( 
\mathbf{0}\right) =\mathbf{I},\ \mathbf{R}_{4}\left( \mathbf{0}\right)
=\left( 
\begin{array}{ccc}
1/\sqrt{2} & 0 & -1/\sqrt{2} \\ 
0 & 1 & 0 \\ 
1/\sqrt{2} & 0 & 1/\sqrt{2}%
\end{array}%
\right) \\
\mathbf{r}_{1}\left( \mathbf{0}\right) &=&\left( -x_{1},0,z_{1}\right) ^{T},%
\mathbf{r}_{2}\left( \mathbf{0}\right) =\left( -x_{2},0,-z_{2}\right) ^{T},%
\mathbf{r}_{3}\left( \mathbf{0}\right) =\left( x_{3},0,z_{3}\right) ^{T} \\
\mathbf{r}_{4}\left( \mathbf{0}\right) &=&\left( x_{4},0,z_{4}\right) ^{T},%
\mathbf{r}_{5}\left( \mathbf{0}\right) =\left( x_{5},0,z_{5}\right) ^{T}.
\end{eqnarray*}%
Therewith the configuration of all bodies are determined by the POE (\ref%
{POEY}). For instance 
\begin{equation*}
\mathbf{C}_{3}\left( \mathbf{q}\right) =\left( 
\begin{array}{cccc}
c_{123} & -s_{123} & 0 & -{d_{2}c_{1}}+({d_{2}}+{d_{3}})c_{12}+({x_{3}}-{%
d_{3}})c_{123} \\ 
s_{123} & c_{123} & 0 & -{d_{2}}s_{1}+({d_{2}}+{d_{3}})s_{12}+({x_{3}}-{d_{3}%
})s_{123} \\ 
0 & 0 & 1 & {z_{3}} \\ 
0 & 0 & 0 & 1%
\end{array}%
\right)
\end{equation*}%
with $c_{123}:=\cos ({q_{1}}+{q_{2}}+{q_{3}})$ etc. The expressions for $%
\mathbf{C}_{4}\left( \mathbf{q}\right) $ and $\mathbf{C}_{5}\left( \mathbf{q}%
\right) $ are rather complicated and are omitted here.

Instead of deducing them from the geometry in figure \ref{figRCMGeom}, the
body-fixed representation of the joint screw coordinates can be determined
with the relation (\ref{XY}). This yields%
\begin{eqnarray*}
{^{1}}\mathbf{X}_{1} &=&\left( 0,0,1,0,-x_{1},0\right) ^{T},{^{2}}\mathbf{X}%
_{2}=\left( 0,0,1,0,d_{2}-x_{2},0\right) ^{T},{^{3}}\mathbf{X}_{3}=\left(
0,0,1,0,-d_{3}+x_{3},0\right) ^{T} \\
{^{4}}\mathbf{X}_{4} &=&(0,0,1,0,-d_{4}/\sqrt{2}-h_{4}/\sqrt{2}+x_{4}/\sqrt{2%
}+z_{4}/\sqrt{2},0)^{T},{^{5}}\mathbf{X}_{5}=\left(
0,0,1,0,-d_{5}+x_{5},0\right) ^{T}.
\end{eqnarray*}%
This example shows the simplicity of the approach.

\section{Velocity of a Kinematic Chain%
\label{secVel}%
}

In this section recursive relations are derived for the four forms of twists
that are introduced in appendix A.3, namely the body-fixed, spatial, hybrid,
and mixed twists \cite{Bruyninckx1996}.

\subsection{Body-fixed Representation of Rigid Body Twists%
\label{secBodyTwist}%
}

\subsubsection{Body-fixed Twists}

\begin{figure}[b]
\centerline{
\includegraphics[width=15cm]{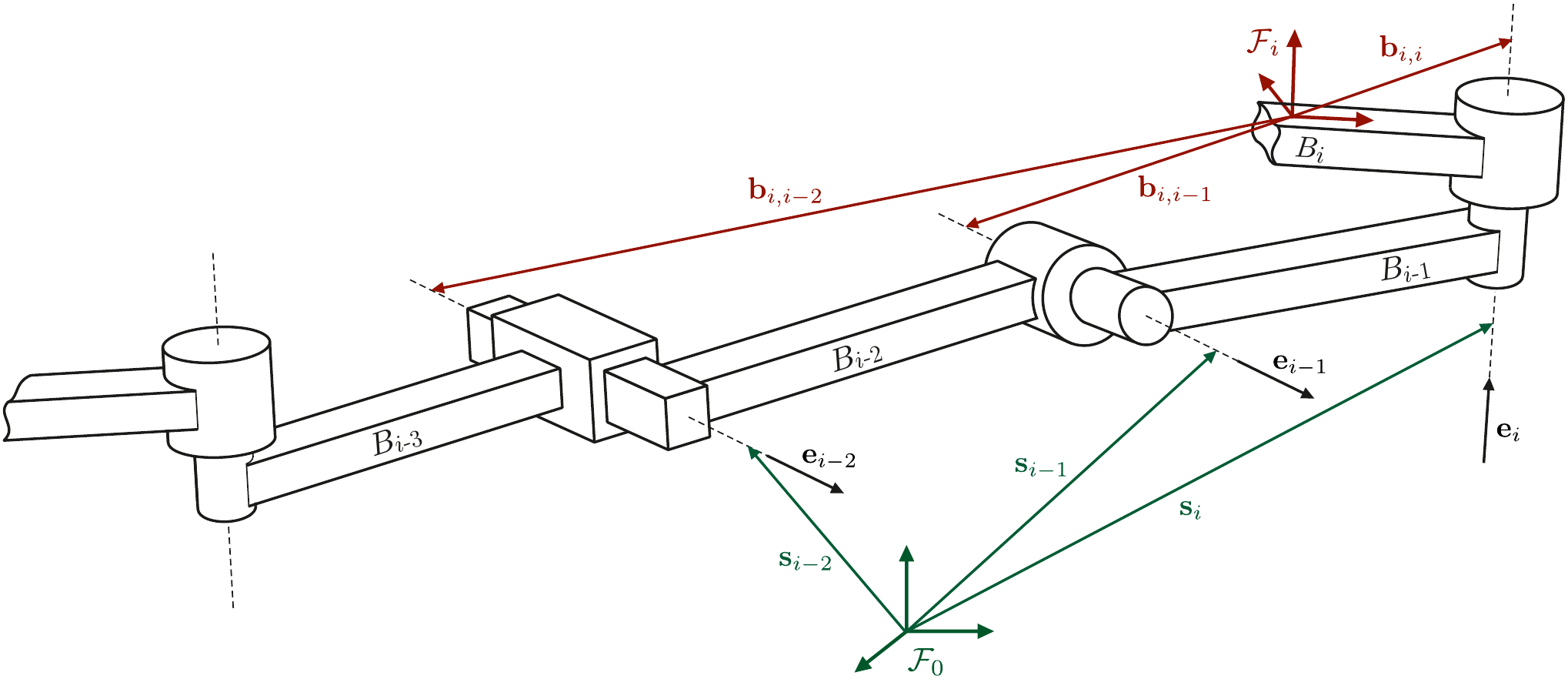}}
\caption{Description of the instantaneous kinematics of a kinematic chain.}
\label{figTwistChain}
\end{figure}
The \emph{body-fixed twist} of body $i$, denoted $\mathbf{V}_{i}^{\text{b}}=(%
\mathbf{\omega }_{i}^{\text{b}},\mathbf{v}_{i}^{\text{b}})^{T}$, is the
aggregate of the angular velocity $\mathbf{\omega }_{i}^{\text{b}}:={^{i}}%
\mathbf{\omega }_{i}$ of its BFR and the translational velocity $\mathbf{v}%
_{i}^{\text{b}}:={^{i}}\mathbf{v}_{i}$ of its origin relative to the IFR
(appendix A.2). The twist of body $i$ in a kinematic chain is the sum of
twists of the joints connecting it to the ground. Represented (measured and
resolved) in its BFR $\mathcal{F}_{i}$ (fig. \ref{figTwistChain}) this is%
\begin{equation}
\mathbf{V}_{i}^{\text{b}}=\dot{q}_{1}\left( 
\begin{array}{c}
{^{i}}\mathbf{e}_{1} \\ 
{^{i}}\mathbf{b}_{i,1}\times {^{i}}\mathbf{e}_{1}+{^{i}}\mathbf{e}_{1}h_{1}%
\end{array}%
\right) +\dot{q}_{2}\left( 
\begin{array}{c}
{^{i}}\mathbf{e}_{2} \\ 
{^{i}}\mathbf{b}_{i,2}\times {^{i}}\mathbf{e}_{2}+{^{i}}\mathbf{e}_{2}h_{2}%
\end{array}%
\right) +\ldots +\dot{q}_{i}\left( 
\begin{array}{c}
{^{i}}\mathbf{e}_{i} \\ 
{^{i}}\mathbf{b}_{i,i}\times {^{i}}\mathbf{e}_{i}+{^{i}}\mathbf{e}_{i}h_{i}%
\end{array}%
\right) .  \label{Vb}
\end{equation}%
Here ${^{i}}\mathbf{b}_{i,j}$ is the instantaneous position vector of a
point on the axis of joint $j$ measured in the BFR $\mathcal{F}_{i}$, and ${%
^{i}}\mathbf{e}_{j}$ is the unit vector along the axis resolved in the BFR.
The instantaneous joint screw coordinates in (\ref{Vb}) are configuration
dependent, and related to the joint screws (\ref{Xb}) and (\ref{XY})
(deduced from reference configuration) by a frame transformation.

\subsubsection{Body-Fixed Jacobian and Recursive Relations}

The body-fixed twist is determined by (\ref{Vbs}) in terms of the
configuration $\mathbf{C}\left( t\right) $. Using (\ref{relConf}) and (\ref%
{derExpX}), from $\mathbf{C}_{i}=\mathbf{C}_{i-1}\mathbf{C}_{i-1,i}$ follows
the recursive relation (notice $\mathbf{C}_{i-1,i}^{-1}=\mathbf{C}_{i,i-1}=%
\mathbf{C}_{i}^{-1}\mathbf{C}_{i-1}$)%
\begin{equation}
\mathbf{V}_{i}^{\text{b}}=\mathbf{Ad}_{\mathbf{C}_{i,i-1}}\mathbf{V}_{i-1}^{%
\text{b}}+{^{i}}\mathbf{X}_{i}\dot{q}_{i}=\mathbf{Ad}_{\mathbf{B}_{i}\exp ({%
^{i}}\mathbf{X}_{i}q_{i})}^{-1}\mathbf{V}_{i-1}^{\text{b}}+{^{i}}\mathbf{X}%
_{i}\dot{q}_{i}.  \label{VbRec}
\end{equation}%
The frame transformations due to the relative motions $\mathbf{C}_{i,i-1}$%
\textbf{\ }of adjacent bodies propagate the twists within the kinematic
chain. The first term on the right hand side of (\ref{VbRec}) is the twist
of body $i-1$ represented in the BFR $\mathcal{F}_{i}$ on body $i$, and the
second term is the additional contribution from joint $i$.

The configuration $\mathbf{C}_{i}$ of body $i$ depends on the joint
variables $q_{j},j\leq i$. The body-fixed twist (\ref{Vbs}) can thus be
expressed as $\widehat{\mathbf{V}}_{i}^{\text{b}}=\sum_{j\leq i}\widehat{%
\mathbf{J}}_{i,j}^{\text{b}}\dot{q}_{j}$ where $\widehat{\mathbf{J}}_{i,j}^{%
\text{b}}:=\mathbf{C}_{i}^{-1}\frac{\partial }{\partial q^{j}}\mathbf{C}_{i}$%
. The POE (\ref{POEX}), together with (\ref{derExpX}) and (\ref{Adexp}),
yields%
\begin{eqnarray}
\mathbf{C}_{i}^{-1}\frac{\partial }{\partial q_{j}}\mathbf{C}_{i} &=&\exp (-{%
^{i}}\widehat{\mathbf{X}}_{i}q_{i})\mathbf{B}_{i}^{-1}\cdots \exp (-{^{j+1}}%
\widehat{\mathbf{X}}_{j+1}q_{j+1})\mathbf{B}_{j+1}^{-1}{^{j}}\widehat{%
\mathbf{X}}_{j}\exp ({^{j+1}}\widehat{\mathbf{X}}_{j+1}q_{j+1})\cdots \exp ({%
^{i}}\widehat{\mathbf{X}}_{i}q_{i})  \notag \\
&=&\mathbf{C}_{i}^{-1}\mathbf{C}_{j}{^{j}}\widehat{\mathbf{X}}_{j}\mathbf{C}%
_{j}^{-1}\mathbf{C}_{i}=\mathrm{Ad}_{\mathbf{C}_{i,j}}({^{j}}\widehat{%
\mathbf{X}}_{j}),\ j\leq i.  \label{derivative1}
\end{eqnarray}%
Using (\ref{XY}) this yields the following relations%
\begin{equation}
\begin{tabular}{|llll|}
\hline
&  &  &  \\ 
& $\mathbf{J}_{i,j}^{\text{b}}%
\hspace{-2ex}%
$ & $=\mathbf{Ad}_{\mathbf{C}_{i,j}}{^{j}}\mathbf{X}_{j}$ &  \\ 
&  &  & $%
\vspace{-2ex}%
$ \\ 
&  & $=\mathbf{Ad}_{\mathbf{C}_{i,j}\mathbf{A}_{j}^{-1}}\mathbf{Y}_{j}$ & 
\\ 
&  &  & $%
\vspace{-2ex}%
$ \\ 
&  & $=\mathbf{Ad}_{\mathbf{C}_{i,j}\mathbf{S}_{j,j}}{^{j-1}}\mathbf{Z}%
_{j},\ \ j\leq i.$ &  \\ 
&  &  &  \\ \hline
\end{tabular}
\label{JbX}
\end{equation}%
The $\mathbf{J}_{i,j}^{\text{b}}$ are the screw coordinate vectors in (\ref%
{Vb}) obtained via a frame transformation (\ref{Ad}) of ${^{j}}\mathbf{X}%
_{j} $ in (\ref{Xb}), or $\mathbf{Y}_{j}$ in (\ref{XY}), to the current
configuration. The body-fixed twist is hence%
\begin{equation}
\mathbf{V}_{i}^{\text{b}}=\sum_{j\leq i}\mathbf{J}_{i,j}^{\text{b}}\dot{q}%
_{j}=\mathsf{J}_{i}^{\text{b}}\dot{\mathbf{q}}.  \label{VbJac}
\end{equation}%
The $6\times n$ matrix%
\begin{equation}
\mathsf{J}_{i}^{\text{b}}=\left( \mathbf{J}_{i,1}^{\text{b}},\ldots ,\mathbf{%
J}_{i,i}^{\text{b}},\mathbf{0}\ldots ,\mathbf{0}\right)  \label{Jbmat}
\end{equation}%
is called the \emph{geometric body-fixed Jacobian} of body $i$ \cite{Murray}%
. It is the central object in all formulations that use body-fixed twists
and Lie group formulations \cite%
{MUBOLieGroup,ParkBobrowPloen1995,Park1994,PloenPark1999}. The geometric
Jacobian appears in the literature under different names. For instance, in 
\cite{maisser1988,maisser1997,MUBOLieGroup} it is called the 'kinematic
basic function (KBF)' as it is the pivotal object for (recursive) evaluation
of MBS kinematics.

The expression (\ref{JbX}) gives rise to the recursive relation%
\begin{equation}
\mathbf{J}_{i,j}^{\text{b}}=\left\{ 
\begin{array}{cl}
\mathbf{Ad}_{\mathbf{C}_{i,i-1}}\mathbf{J}_{i-1,j}^{\text{b}}, & j<i%
\vspace{-1ex}
\\ 
&  \\ 
{^{i}}\mathbf{X}_{i} & j=i\ \ .%
\end{array}%
\right.  \label{JbRec}
\end{equation}%
This is essentially another form of the recursion (\ref{VbRec}), using (\ref%
{VbJac}).

\begin{remark}[Dependence on joint variables]
With (\ref{POEX}), respectively (\ref{POEY}), it is clear that the Jacobian $%
\mathbf{J}_{i}^{\text{b}}$ of body $i$ can only depend on $q_{1},\ldots
,q_{i}$. Moreover, noting in (\ref{JbX}) that $\mathbf{C}_{i,j}=\mathbf{C}%
_{i}^{-1}\mathbf{C}_{j}$ is independent from $q_{1}$, it follows that $%
\mathbf{J}_{i}^{\text{b}}$ depends on $q_{2},\ldots ,q_{i}$, i.e. it is
independent from the first joint in the chain. This is obvious from a
kinematic perspective since $\mathbf{V}_{i}^{\text{b}}$ is the sum of twists
of the preceding bodies in the chain expressed in the BFR on body $i$. This
only depends on the configuration of the bodies relative to body $i$ but not
on the absolute configuration of the overall chain, which is determined by $%
q_{1}$.
\end{remark}

\begin{remark}[Required data]
The second form in (\ref{JbX}) in conjunction with (\ref{POEY}) allows for
computation of the body-fixed Jacobian without introducing body-fixed JFRs.
The only information needed are the joint screw coordinates $\mathbf{Y}_{j}$
represented in the IFR and the reference configurations $\mathbf{A}_{j}$.
\end{remark}

\begin{remark}[Change of reference frame]
When another BFR on body $i$ is used, which is related to the original BFR
by $\mathbf{S}\in SE\left( 3\right) $, its configuration is given by $%
\mathbf{C}_{i}^{\prime }=\mathbf{C}_{i}\mathbf{S}$. The corresponding
body-fixed twist follows from (\ref{Vbs}) as $\widehat{\mathbf{V}}%
_{i}^{\prime \text{b}}=\mathbf{C}_{i}^{\prime -1}\dot{\mathbf{C}}%
_{i}^{\prime }=\mathbf{S}^{-1}\mathbf{C}_{i}^{-1}\dot{\mathbf{C}}_{i}\mathbf{%
S}=\mathrm{Ad}_{\mathbf{S}}^{-1}(\widehat{\mathbf{V}}_{i}^{\text{b}})$, and
in vector form%
\begin{equation}
\mathbf{V}_{i}^{\prime \text{b}}=\mathbf{Ad}_{\mathbf{S}}^{-1}\mathbf{V}%
_{i}^{\text{b}}.
\end{equation}%
On the other hand the body-fixed twist is invariant under a change of IFR,
which is given by $\mathbf{C}_{i}^{\prime }=\mathbf{SC}_{i}$. Body-fixed
twists are therefore called left-invariant vector fields on $SE\left(
3\right) $ since left multiplication of $\mathbf{C}_{i}$ with any $\mathbf{S}%
\in SE\left( 3\right) $ does not affect $\mathbf{V}_{i}^{\text{b}}$.
\end{remark}

\begin{remark}[Application of body-fixed representation]
The recursive relations for body-fixed twist and Jacobian are the basis for
the MBS dynamics algorithms in \cite%
{Anderson1992,AndersonCritchley2003,Bae2001,Fijany1995,Hollerbach1980,LillyOrin1991,Liu1988,ParkBobrowPloen1995,Park1994,PloenPark1999,VukobratovicPotkonjak1979}%
. In \cite{Anderson1992,AndersonCritchley2003} the adjoint transformation
matrix in (\ref{VbRec}) was called the 'shift matrix', and $\mathbf{X}_{i}$
was called the 'motion map matrix'. However, the geometric background was
rarely exploited as in \cite{ParkBobrowPloen1995,Park1994,PloenPark1999} and 
\cite{Liu1988}. Remarkably Liu \cite{Liu1988} already presented all relevant
formulations in terms of screws.
\end{remark}

\subsubsection{Body-Fixed System Jacobian and its Decomposition%
\label{secFactorB}%
}

The body-fixed twists are summarized in the overall twist vector $\mathsf{V}%
^{\text{b}}=(\mathbf{V}_{1}^{\text{b}},\ldots ,\mathbf{V}_{n}^{\text{b}%
})^{T} $. The recursion (\ref{VbRec}) can then be written in matrix form%
\begin{equation}
\mathsf{V}^{\text{b}}=\mathsf{D}^{\text{b}}\mathsf{V}^{\text{b}}+\mathsf{X}^{%
\text{b}}\dot{\mathbf{q}}  \label{Vsb1}
\end{equation}%
with%
\begin{equation}
\mathsf{D}^{\text{b}}:=\left( 
\begin{array}{ccccc}
\mathbf{0} & \mathbf{0} & \mathbf{0} &  & \mathbf{0} \\ 
\mathbf{Ad}_{\mathbf{C}_{2,1}} & \mathbf{0} & \mathbf{0} & \cdots &  \\ 
\mathbf{0} & \mathbf{Ad}_{\mathbf{C}_{3.2}} & \mathbf{0} &  &  \\ 
\vdots & \vdots & \ddots & \ddots &  \\ 
\mathbf{0} & \mathbf{0} & \cdots & \mathbf{Ad}_{\mathbf{C}_{n,n-1}} & 
\mathbf{0}%
\end{array}%
\right) ,\ \mathsf{X}^{\text{b}}:=\mathrm{diag}~\left( {^{1}\mathbf{X}}%
_{1},\ldots ,{^{n}}\mathbf{X}_{n}\right) .  \label{Db}
\end{equation}%
On the other hand, the recursive expression for the Jacobian (\ref{JbRec})
reads in matrix form%
\begin{equation}
\mathsf{V}^{\text{b}}=\mathsf{J}^{\text{b}}\dot{\mathbf{q}}=\mathsf{A}^{%
\text{b}}\mathsf{X}^{\text{b}}\dot{\mathbf{q}}  \label{VbAX}
\end{equation}%
where the $6n\times n$ matrix $\mathsf{J}^{\text{b}}=\mathsf{A}^{\text{b}}%
\mathsf{X}^{\text{b}}$ is the \emph{system Jacobian in body-fixed
representation}, and 
\begin{equation}
\mathsf{A}^{\text{b}}:=\left( 
\begin{array}{ccccc}
\mathbf{I} & \mathbf{0} & \mathbf{0} &  & \mathbf{0} \\ 
\mathbf{Ad}_{\mathbf{C}_{2,1}} & \mathbf{I} & \mathbf{0} & \cdots & \mathbf{0%
} \\ 
\mathbf{Ad}_{\mathbf{C}_{3,1}} & \mathbf{Ad}_{\mathbf{C}_{3,2}} & \mathbf{I}
&  & \mathbf{0} \\ 
\vdots & \vdots & \ddots & \ddots &  \\ 
\mathbf{Ad}_{\mathbf{C}_{n,1}} & \mathbf{Ad}_{\mathbf{C}_{n,2}} & \cdots & 
\mathbf{Ad}_{\mathbf{C}_{n,n-1}} & \mathbf{I}%
\end{array}%
\right)  \label{Ab}
\end{equation}%
is the screw transformation matrix. Comparing (\ref{Vsb1}) and (\ref{VbAX})
shows that $\mathsf{A}^{\text{b}}=(\mathbf{I}-\mathsf{D}^{\text{b}})^{-1}$.
In fact $\mathsf{D}^{\text{b}}$ is nilpotent so that the von-Neumann series%
\begin{equation}
\mathsf{A}^{\text{b}}=(\mathbf{I}-\mathsf{D}^{\text{b}})^{-1}=\mathbf{I}+%
\mathsf{D}^{\text{b}}+(\mathsf{D}^{\text{b}})^{2}+\ldots +(\mathsf{D}^{\text{%
b}})^{n}  \label{Invsb}
\end{equation}%
terminates with $(\mathsf{D}^{\text{b}})^{n+1}=0$. That is, $\mathbf{A}^{%
\text{b}}$ is the 1-resolvent of $\mathsf{D}^{\text{b}}$, which is the
fundamental point of departure for many $O\left( n\right) $ algorithms.
Moreover, (\ref{Invsb}) is another form of the recursive coordinate
transformations. Hence the inverse $(\mathsf{A}^{\text{b}})^{-1}=(\mathbf{I}-%
\mathsf{D}^{\text{b}})$.

\begin{remark}[Overall inverse kinematics solution]
\label{remInvKin}%
The above result allows for a simple transformation from body-fixed
velocities to the corresponding joint rates. When the twists of all bodies
are given, (\ref{VbAX}) is an overdetermined linear system in $\dot{\mathbf{q%
}}$. It has a unique solution as long as the twists are consistent with the
kinematics. Premultiplication of (\ref{VbAX}) with $(\mathbf{I}-\mathsf{D}^{%
\text{b}})$ followed by $(\mathsf{X}^{\text{b}})^{T}$ and $((\mathsf{X}^{%
\text{b}})^{T}\mathsf{X}^{\text{b}})^{-1}$ yields%
\begin{equation}
\dot{\mathbf{q}}=((\mathsf{X}^{\text{b}})^{T}\mathsf{X}^{\text{b}})^{-1}(%
\mathsf{X}^{\text{b}})^{T}(\mathbf{I}-\mathsf{D}^{\text{b}})\mathsf{V}^{%
\text{b}}.  \label{solqdot1}
\end{equation}%
The $n\times n$ diagonal matrix $(\mathsf{X}^{\text{b}T}\mathsf{X}^{\text{b}%
})^{-1}=\mathrm{diag}(1/\left\Vert {^{1}\mathbf{X}}_{1}\right\Vert
^{2},\ldots ,1/\left\Vert {^{n}\mathbf{X}}_{n}\right\Vert ^{2})$ has full
rank. Due to the block diagonal structure this yields the solutions $\dot{q}%
_{i}={^{i}\mathbf{X}}_{i}^{T}(\mathbf{V}_{i}^{\text{b}}-\mathbf{Ad}_{\mathbf{%
C}_{i,i-1}}\mathbf{V}_{i-1}^{\text{b}})/\left\Vert {^{i}\mathbf{X}}%
_{i}\right\Vert ^{2}$ for the individual joints. This is indeed the
projection of the relative twist of body $i$ w.r.t. body $i-1$ onto the axis
of joint $i$. It is an exact solution of the inverse kinematics for the
overall MBS, presumed that the twists are compatible, i.e. satisfy (\ref%
{VbRec}). If this is not the case, (\ref{solqdot1}) is the unique
pseudoinverse solution of the system (\ref{VbAX}) of $6n$ equations for the $%
n$ unknown $\dot{q}_{i}$ minimizing the residual error. This can be
considered as the \emph{generalized inverse kinematics problem}: given
desired twists of all individual links, find the joint rates that best
reproduce these twists. This can be applied, for instance, to the inverse
kinematics of human body models processing motion capture data (estimated
position and orientation of body segments) and when noisy data is processed.

While the solution (\ref{solqdot1}) seems straightforward, it should be
remarked that there is no frame invariant inner product on $se\left(
3\right) $, i.e. no norm of screws can be defined that is invariant under a
change of reference frame \cite{Selig}. The correctness of (\ref{solqdot1})
follows by regarding the transposed joint screw coordinates as co-screws,
and ${^{i}\mathbf{X}}_{i}^{T}{^{i}\mathbf{X}}_{i}$ is the pairing of screw
and co-screw coordinates rather than an inner product.
\end{remark}

\subsection{Spatial Representation of Rigid Body Twists%
\label{secSpatialTwists}%
}

\subsubsection{Spatial Twists}

A representation of the body twist, which is less common in MBS modeling but
frequently used in mechanism theory, is the so-called \emph{spatial twist}
denoted $\mathbf{V}_{i}^{\text{s}}=(\mathbf{\omega }_{i}^{\text{s}},\mathbf{v%
}_{i}^{\text{s}})^{T}$. This is the twist of body $i$ represented in the
IFR. It consists of the angular velocity of the BFR of body $i$ measured and
resolved in the IFR, and the translational velocity $\mathbf{v}_{i}^{\text{s}%
}:=\dot{\mathbf{r}}_{i}-\mathbf{\omega }_{i}^{\text{s}}\times \mathbf{r}_{i}$
of the (possibly imaginary) point on the body that is momentarily traveling
through the origin of the IFR measured and resolved in the IFR (appendix
A.2). With the notation in fig. \ref{figTwistChain}, the spatial twist of
body $i$ is geometrically readily constructed as%
\begin{equation}
\mathbf{V}_{i}^{\text{s}}=\dot{q}_{1}\left( 
\begin{array}{c}
\mathbf{e}_{1} \\ 
\mathbf{s}_{1}\times \mathbf{e}_{1}+h_{1}\mathbf{e}_{1}%
\end{array}%
\right) +\dot{q}_{2}\left( 
\begin{array}{c}
\mathbf{e}_{2} \\ 
\mathbf{s}_{2}\times \mathbf{e}_{2}+h_{2}\mathbf{e}_{2}%
\end{array}%
\right) +\ldots +\dot{q}_{i}\left( 
\begin{array}{c}
\mathbf{e}_{i} \\ 
\mathbf{s}_{i}\times \mathbf{e}_{i}+h_{i}\mathbf{e}_{i}%
\end{array}%
\right)  \label{Vsi}
\end{equation}%
where $\mathbf{s}_{j}$ is the position vector of a point on the joint axis $%
j $ expressed in the IFR. The screw coordinates in (\ref{Vsi}) are
configuration dependent. They are equal to $\mathbf{Y}_{j}$ in the reference
configuration $\mathbf{q}=\mathbf{0}$, where $\mathbf{s}_{i}=\mathbf{y}_{i}$.

\subsubsection{Spatial Jacobian and Recursive Relations}

In order to derive an analytic expression, using the POE, the definition (%
\ref{Vbs}) of the spatial twist is applied. As apparent from (\ref{Vsi}),
the non-vanishing instantaneous joint screws are identical for all bodies.
This is clear since the IFR is the only reference frame involved. The
spatial twist can thus be expressed as $\mathbf{V}_{i}^{\text{s}%
}=\sum_{j\leq i}\mathbf{J}_{j}^{\text{s}}\dot{q}_{j}$ with $\widehat{\mathbf{%
J}}_{j}^{\text{s}}:=\frac{\partial }{\partial q_{j}}\mathbf{C}_{i}\mathbf{C}%
_{i}^{-1}$. Using the POE, a straightforward derivation analogous to (\ref%
{derivative1}) yields%
\begin{equation}
\begin{tabular}{|llll|}
\hline
&  &  &  \\ 
& $\mathbf{J}_{j}^{\text{s}}%
\hspace{-2ex}%
$ & $=\mathbf{Ad}_{\mathbf{C}_{j}}{^{j}}\mathbf{X}_{j}$ &  \\ 
&  &  & $%
\vspace{-2ex}%
$ \\ 
&  & $=\mathbf{Ad}_{\mathbf{C}_{j}\mathbf{A}_{j}^{-1}}\mathbf{Y}_{j}$ &  \\ 
&  &  & $%
\vspace{-2ex}%
$ \\ 
&  & $=\mathbf{Ad}_{\mathbf{C}_{j}\mathbf{S}_{j,j}}{^{j-1}}\mathbf{Z}_{j},\
\ j\leq i.$ &  \\ 
&  &  &  \\ \hline
\end{tabular}
\label{JsX}
\end{equation}%
The $\mathbf{J}_{j}^{\text{s}}$ is the instantaneous screw coordinate vector
of joint $j$ in (\ref{Vsi}) in spatial representation, i.e. represented in
the IFR. The matrix%
\begin{equation}
\mathsf{J}_{i}^{\text{s}}=%
\Big%
(\mathbf{J}_{1}^{\text{s}},\ldots ,\mathbf{J}_{i}^{\text{s}},\mathbf{0}%
,\ldots ,\mathbf{0}%
\Big%
)  \label{Js}
\end{equation}%
is called the \emph{spatial Jacobian}\textit{\ }of body $i$. The relations (%
\ref{Vsi}) and (\ref{JsX}) yield the following recursive expression for the
spatial twists of bodies in a kinematic chain%
\begin{equation}
\mathbf{V}_{i}^{\text{s}}=\mathbf{V}_{i-1}^{\text{s}}+\mathbf{J}{_{i}^{\text{%
s}}}\dot{q}_{i}.  \label{Vsrec}
\end{equation}

\begin{remark}
\label{remJs}%
The spatial representation has remarkable advantages. The velocity recursion
(\ref{Vsrec}) is the simplest possible since the twists of individual bodies
can simply be added \emph{without any coordinate transformation}. An
important observation is that $\mathbf{J}_{j}^{\text{s}}$ is intrinsic to
the joint $j$. The non-zero screw vectors in the Jacobian (\ref{Js}) are
thus the \emph{same for all bodies}. This is a consequence of using a single
spatial reference frame.
\end{remark}

\subsubsection{Spatial System Jacobian and its Decomposition}

The overall spatial twist $\mathsf{V}^{\text{s}}=(\mathbf{V}_{1}^{\text{s}%
},\ldots ,\mathbf{V}_{n}^{\text{s}})^{T}$ of the kinematic chain is
determined as%
\begin{equation}
\mathsf{V}^{\text{s}}=\mathsf{J}^{\text{s}}\dot{\mathbf{q}}
\end{equation}%
where the \emph{spatial system Jacobian} possesses the factorizations%
\begin{equation}
\mathsf{J}^{\text{s}}=\mathsf{A}^{\text{s}}\mathsf{Y}^{\text{s}}=\mathsf{A}^{%
\text{sb}}\mathsf{X}^{\text{b}}=\mathsf{A}^{\text{sh}}\mathsf{X}^{\text{h}}.
\label{VsJac}
\end{equation}%
Therein it is $\mathsf{Y}^{\text{s}}=\mathrm{diag}~\left( {\mathbf{Y}}%
_{1},\ldots ,\mathbf{Y}_{n}\right) $, $\mathsf{X}^{\text{h}}$ in (\ref{Asb}%
), and%
\begin{eqnarray}
\mathsf{A}^{\text{sb}}:= &&\left( 
\begin{array}{ccccc}
\mathbf{Ad}_{\mathbf{C}_{1}} & \mathbf{0} & \mathbf{0} &  & \mathbf{0} \\ 
\mathbf{Ad}_{\mathbf{C}_{1}} & \mathbf{Ad}_{\mathbf{C}_{2}} & \mathbf{0} & 
& \mathbf{0} \\ 
\mathbf{Ad}_{\mathbf{C}_{1}} & \mathbf{Ad}_{\mathbf{C}_{2}} & \mathbf{Ad}_{%
\mathbf{C}_{3}} &  & \mathbf{0} \\ 
\vdots  & \vdots  & \ddots  & \ddots  & \vdots  \\ 
\mathbf{Ad}_{\mathbf{C}_{1}} & \mathbf{Ad}_{\mathbf{C}_{2}} & \cdots  & 
\mathbf{Ad}_{\mathbf{C}_{n-1}} & \mathbf{Ad}_{\mathbf{C}_{n}}%
\end{array}%
\right) ,\ \ \ \mathsf{A}^{\text{sh}}:=\left( 
\begin{array}{ccccc}
\mathbf{Ad}_{\mathbf{r}_{1}} & \mathbf{0} & \mathbf{0} &  & \mathbf{0} \\ 
\mathbf{Ad}_{\mathbf{r}_{1}} & \mathbf{Ad}_{\mathbf{r}_{2}} & \mathbf{0} & 
& \mathbf{0} \\ 
\mathbf{Ad}_{\mathbf{r}_{1}} & \mathbf{Ad}_{\mathbf{r}_{2}} & \mathbf{0} & 
& \mathbf{0} \\ 
\vdots  & \vdots  & \ddots  & \ddots  &  \\ 
\mathbf{Ad}_{\mathbf{r}_{1}} & \mathbf{Ad}_{\mathbf{r}_{2}} & \cdots  & 
\mathbf{Ad}_{\mathbf{r}_{n-1}} & \mathbf{Ad}_{\mathbf{r}_{n}}%
\end{array}%
\right) \ \ \ \   \label{As} \\
\mathsf{A}^{\text{s}}:= &&\mathsf{A}^{\text{sb}}\mathrm{diag}~(\mathbf{Ad}_{%
\mathbf{A}_{1}}^{-1},\ldots ,\mathbf{Ad}_{\mathbf{A}_{n}}^{-1}).  \notag
\end{eqnarray}%
All non-zero entries in a column of these matrices are identical. Hence the
construction of these matrices only requires determination of the $n$
entries in the last row that are copied into the upper triangular block
matrix. The factorization (\ref{VsJac}) gives rise to an expression for its
inverse. Noting $\mathsf{A}^{\text{sb}}=\mathrm{diag}~(\mathbf{Ad}_{\mathbf{C%
}_{1}},\ldots ,\mathbf{Ad}_{\mathbf{C}_{n}})\mathsf{A}^{\text{b}}$ the
relation for the inverse of $\mathsf{A}^{\text{b}}$ in terms of matrix $%
\mathsf{D}^{\text{b}}$ in (\ref{Db}) yields%
\begin{eqnarray}
(\mathsf{A}^{\text{sb}})^{-1} &=&(\mathbf{I}-\mathsf{D}^{\text{b}})\mathrm{%
diag}~(\mathbf{Ad}_{\mathbf{C}_{1}},\ldots ,\mathbf{Ad}_{\mathbf{C}_{n}}) 
\notag \\
&=&\left( 
\begin{array}{ccccc}
\mathbf{Ad}_{\mathbf{C}_{1}}^{-1} & \mathbf{0} & \mathbf{0} &  & \mathbf{0}
\\ 
-\mathbf{Ad}_{\mathbf{C}_{2}} & \mathbf{Ad}_{\mathbf{C}_{2}}^{-1} & \mathbf{0%
} &  & \mathbf{0} \\ 
\mathbf{0} & -\mathbf{Ad}_{\mathbf{C}_{3}} & \mathbf{Ad}_{\mathbf{C}%
_{3}}^{-1} &  & \mathbf{0} \\ 
\vdots  & \vdots  & \ddots  & \ddots  & \vdots  \\ 
\mathbf{0} & \mathbf{0} & \cdots  & -\mathbf{Ad}_{\mathbf{C}_{n}} & \mathbf{%
Ad}_{\mathbf{C}_{n}}^{-1}%
\end{array}%
\right)   \label{Asbinv}
\end{eqnarray}%
and $(\mathsf{A}^{\text{s}})^{-1}$ accordingly.

\begin{remark}[Dependence on joint variables]
Similarly to the body-fixed twist, it follows from $\mathbf{J}_{i}^{\text{s}%
}=\mathbf{Ad}_{\mathbf{C}_{i}}{^{i}}\mathbf{X}_{i}=\mathbf{Ad}_{\mathbf{C}%
_{i-1}}\mathbf{Ad}_{\mathbf{B}_{i}\exp {^{i}}\mathbf{X}_{i}q_{i}}{^{i}}%
\mathbf{X}_{i}=\mathbf{Ad}_{\mathbf{C}_{i-1}}{^{i}}\mathbf{X}_{i}$ being
independent from $q_{i}$, that the spatial Jacobian of body $i$ only depends
on $q_{1},\ldots ,q_{i-1}$. Indeed, the motion of joint $i$ does not change
its screw axis about which body $i$ is moving.
\end{remark}

\begin{remark}[Change of reference frame]
The spatial twist is called a right-invariant vector field on $SE\left(
3\right) $ because it does not change when $\mathbf{C}_{i}$ is
postmultiplied by any $\mathbf{S}\in SE\left( 3\right) $, representing a
change of body-fixed RFR. Under a change of IFR according to $\mathbf{C}%
_{i}^{\prime }=\mathbf{SC}_{i}$ the spatial twists transform as%
\begin{equation}
\mathbf{V}_{i}^{\prime \text{s}}=\mathbf{Ad}_{\mathbf{S}}\mathbf{V}_{i}^{%
\text{s}}.
\end{equation}
\end{remark}

\begin{remark}[Application of spatial representation]
The spatial twist is used almost exclusively in mechanism kinematics (often
without mentioning it) but is becoming accepted for MBS modeling since it
was introduced in \cite{Featherstone1983,Featherstone2008}. For kinematic
analysis of mechanisms it is common practice to (instantaneously) locate the
global reference frame so that it coincides with the frame where
kinetostatic properties (twists, wrenches) are observed, usually at the
end-effector. For a serial robotic manipulator the end-effector frame is
located at the terminal link of the chain, so that $\mathbf{A}_{n}=\mathbf{I}
$, and $\mathbf{V}_{n}^{\text{s}}$ is then the spatial end-effector twist.
From their definition follows that the spatial and hybrid twist (see next
section) of body $i$ are numerically identical when the BFR $\mathcal{F}_{i}$
overlaps with the IFR $\mathcal{F}_{0}$.\newline
The most prominent use of the spatial representation in dynamics is the $%
O\left( n\right) $ forward dynamics method by Featherstone \cite%
{Featherstone1983,Featherstone2008}. This has not yet been widely applied in
MBS dynamics. This may be due to use of an uncommon choice of reference
point (the IFR origin) at which the spatial entities are measured, so that
results and interaction wrenches must be transformed to body-fixed reference
frames. The spatial representation of twists must not be confused with the
'spatial vector' notation proposed in \cite%
{Featherstone1983,Featherstone2008}. The latter is a general expression of
twists as 6-vectors (like body-fixed and spatial) but without reference to a
particular frame in which the components are resolved. This allows for
abstract derivation of kinematic relations, but these relation must
eventually be resolved in a particular frame, and this eventually determines
the computational effort.\newline
A notable application of the spatial twist is the modeling and numerical
integration of non-linear elastic MBS where it is called the base pole
velocity \cite{Borri2001a} or fixed pole velocity \cite{Bottasso1998} and
the intrinsic coupling of translational and angular velocity (according to
the screw motion) was discussed. The corresponding momentum balance and
conservation properties are discussed in \cite{Borri2001b,Borri2003} (see
also \cite{Part2}).
\end{remark}

\begin{remark}
As in remark \ref{remInvKin}, the relation (\ref{Asbinv}) gives rise to an
overall inverse kinematics solution. For given spatial twists $\mathbf{V}%
_{i}^{\text{s}}$ this reads in components $\dot{q}_{i}={^{i}\mathbf{X}}%
_{i}^{T}\mathbf{Ad}_{\mathbf{C}_{i}}^{-1}(\mathbf{V}_{i}^{\text{s}}-\mathbf{V%
}_{i-1}^{\text{s}})/\left\Vert {^{i}\mathbf{X}}_{i}\right\Vert ^{2}$.
\end{remark}

\subsection{Hybrid Form of Rigid Body Twists}

\subsubsection{Hybrid Twists}

In various applications it is beneficial to measure the twist of a body in
the body-fixed BFR but resolve it in the IFR. This is commonly referred to
as the \emph{hybrid twist} \cite{Bruyninckx1996,Murray}, denoted $\mathbf{V}%
_{i}^{\text{h}}=(\mathbf{\omega }_{i}^{\text{s}},\dot{\mathbf{r}}_{i})^{T}$.
The geometric construction (fig. \ref{figTwistChain}) yields%
\begin{equation}
\mathbf{V}_{i}^{\text{h}}=\dot{q}_{1}\left( 
\begin{array}{c}
\mathbf{e}_{1} \\ 
\mathbf{b}_{i,1}\times \mathbf{e}_{1}+h_{1}\mathbf{e}_{1}%
\end{array}%
\right) +\dot{q}_{2}\left( 
\begin{array}{c}
\mathbf{e}_{2} \\ 
\mathbf{b}_{i,2}\times \mathbf{e}_{2}+h_{2}\mathbf{e}_{2}%
\end{array}%
\right) +\ldots +\dot{q}_{i}\left( 
\begin{array}{c}
\mathbf{e}_{i} \\ 
\mathbf{b}_{i,i}\times \mathbf{e}_{i}+h_{i}\mathbf{e}_{i}%
\end{array}%
\right) .  \label{Vbsi}
\end{equation}%
As in (\ref{Vb}), $\mathbf{b}_{i,j}$ is the position vector of a point on
the axis of joint $j$ measured from the BFR $\mathcal{F}_{i}$ of body $i$,
and $\mathbf{e}_{j}$ is a unit vector along the axis, but now expressed in
the IFR $\mathcal{F}_{0}$. This was originally introduced in \cite%
{Waldron1982} and \cite{Whitney1972} and is used in various $O\left(
n\right) $ dynamics algorithms (remark \ref{remAppHyb}).

\subsubsection{Hybrid Jacobian and Recursive Relations}

The hybrid twist is merely the body-fixed twist resolved in the IFR. Using (%
\ref{AdRr}) this transformation is $\mathbf{V}_{i}^{\text{h}}=\mathbf{Ad}_{%
\mathbf{R}_{i}}\mathbf{V}_{i}^{\text{b}}$ where $\mathbf{R}_{i}$ is the
rotation matrix of body $i$. Then (\ref{VbJac}) leads to%
\begin{equation}
\mathbf{V}_{i}^{\text{h}}=\sum_{j\leq i}\mathbf{J}_{i,j}^{\text{h}}\dot{q}%
_{j}=\mathsf{J}_{i}^{\text{h}}\dot{\mathbf{q}}
\end{equation}%
with the columns $\mathbf{J}_{i,j}^{\text{h}}:=\mathbf{Ad}_{\mathbf{R}_{i}}%
\mathbf{J}_{i,j}^{\text{b}}$ of the \emph{hybrid Jacobian} $\mathsf{J}_{i}^{%
\text{h}}=\left( \mathbf{J}_{i,1}^{\text{h}},\ldots ,\mathbf{J}_{i,i}^{\text{%
h}},\mathbf{0}\ldots ,\mathbf{0}\right) $. The recursive expressions (\ref%
{JbRec}) and (\ref{VbRec}) remain valid when all screw coordinate vectors
are resolved in the IFR. The joint screw coordinates are then configuration
dependent. The screw coordinate vector of joint $j$ measured in the BFR $%
\mathcal{F}_{j}$ on body $j$ and resolved in the IFR $\mathcal{F}_{0}$ is
related to (\ref{Xb}) via%
\begin{equation}
{^{0}}\mathbf{X}_{j}^{j}=\mathbf{Ad}_{\mathbf{R}_{j}}{^{j}}\mathbf{X}%
_{j}=\left( 
\begin{array}{c}
\mathbf{e}_{j} \\ 
\mathbf{x}_{j,j}\times \mathbf{e}_{j}+h_{j}\mathbf{e}_{j}%
\end{array}%
\right) .  \label{Xh}
\end{equation}%
As in (\ref{Xb}), the position vector $\mathbf{x}_{j,j}$ of a point on the
axis of joint $j$ measured from the BFR $\mathcal{F}_{j}$ but now resolved
in the IFR. The relations $\mathbf{Ad}_{\mathbf{C}_{i}}=\mathbf{Ad}_{\mathbf{%
r}_{i}}\mathbf{Ad}_{\mathbf{R}_{i}}$ and $\mathbf{r}_{i,j}=\mathbf{r}_{j}-%
\mathbf{r}_{i}$ lead to%
\begin{eqnarray}
\mathbf{Ad}_{\mathbf{R}_{i}}\mathbf{J}_{i,j}^{\text{b}} &=&\mathbf{Ad}_{%
\mathbf{R}_{i}}\mathbf{Ad}_{\mathbf{C}_{i,j}}{^{j}}\mathbf{X}_{j}=\mathbf{Ad}%
_{\mathbf{R}_{i}}\mathbf{Ad}_{\mathbf{C}_{i,j}}\mathbf{Ad}_{\mathbf{R}%
_{j}}^{-1}{^{0}}\mathbf{X}_{j}^{j}=\mathbf{Ad}_{\mathbf{R}_{i}}\mathbf{Ad}_{%
\mathbf{C}_{i}}^{-1}\mathbf{Ad}_{\mathbf{C}_{j}}\mathbf{Ad}_{\mathbf{R}%
_{j}}^{-1}{^{0}}\mathbf{X}_{j}^{j}  \notag \\
&=&\mathbf{Ad}_{\mathbf{r}_{i}}^{-1}\mathbf{Ad}_{\mathbf{r}_{j}}{^{0}}%
\mathbf{X}_{j}^{j}=\mathbf{Ad}_{-\mathbf{r}_{i}}\mathbf{Ad}_{\mathbf{r}_{j}}{%
^{0}}\mathbf{X}_{j}^{j}=\mathbf{Ad}_{\mathbf{r}_{j}-\mathbf{r}_{i}}{^{0}}%
\mathbf{X}_{j}^{j}=\mathbf{Ad}_{\mathbf{r}_{i,j}}{^{0}}\mathbf{X}_{j}^{j}.
\label{temp2}
\end{eqnarray}%
Therewith the columns of the hybrid Jacobian of body $i$ are%
\begin{equation}
\begin{tabular}{|lll|}
\hline
&  &  \\ 
& $\mathbf{J}_{i,j}^{\text{h}}=\mathbf{Ad}_{\mathbf{r}_{i,j}}{^{0}}\mathbf{X}%
_{j}^{j},\ \ j\leq i.$ &  \\ 
&  &  \\ \hline
\end{tabular}
\label{Jh}
\end{equation}%
The $\mathbf{J}_{i,j}^{\text{h}}$ is the instantaneous screw coordinate of
joint $j$ in (\ref{Vbsi}) measured at BFR on body $i$ and resolved in the
IFR. In the hybrid form, all vectors are resolved in the IFR. That is, the
screw coordinates ${^{0}}\mathbf{X}_{j}^{j}$ depend on the current
configuration $\mathbf{q}$ even though the joint axis is constant within
body $j$. The hybrid twist is resolved in the IFR. Since the screw
coordinates ${^{0}}\mathbf{X}_{j}^{j}$ are already resolved in the IFR the
transformation to the current configuration, in order to determine the
instantaneous joint screws $\mathbf{J}_{i,j}^{\text{h}}\left( \mathbf{q}%
\right) $, only requires translations of origins. This is obtained by
shifting the reference point according to $\mathbf{r}_{i,j}$, which is why
the matrix $\mathbf{Ad}_{\mathbf{r}_{i,j}}$ is also called the 'shift dyad' 
\cite{CritchleyAnderson2003}. This is not a frame transformation.

The relation (\ref{Jh}) gives rise to the recursive relation for the hybrid
Jacobian%
\begin{equation}
\mathbf{J}_{i,j}^{\text{h}}=\left\{ 
\begin{array}{cl}
\mathbf{Ad}_{\mathbf{r}_{i,i-1}}\mathbf{J}_{i-1,j}^{\text{h}}, & j<i%
\vspace{-1ex}
\\ 
&  \\ 
{^{0}}\mathbf{X}_{i}^{i}, & j=i\ .%
\end{array}%
\right. \ \   \label{recJh2}
\end{equation}%
and, analogous to (\ref{VbRec}), for the hybrid twists%
\begin{equation}
\mathbf{V}_{i}^{\text{h}}=\mathbf{Ad}_{\mathbf{r}_{i,i-1}}\mathbf{V}_{i-1}^{%
\text{h}}+{^{0}}\mathbf{X}_{i}^{i}\dot{q}^{i}.  \label{VhRec}
\end{equation}%
The advantage of the hybrid form over the body-fixed is that (\ref{recJh2})
only involves the relative displacement $\mathbf{r}_{i,i-1}$ in contrast to
the complete relative configuration $\mathbf{C}_{i,i-1}$ in (\ref{JbRec}).
It must be recalled, however, that the vectors $\mathbf{e}_{j}$ and $\mathbf{%
r}_{j,j}$ must be transformed to the IFR. Furthermore, when formulating
equations of motion, the inertia properties of the body must be resolved in
the IFR so that they become configuration dependent \cite{Part2}.

\begin{remark}[Application of hybrid representation]
\label{remAppHyb}%
The hybrid form was used in \cite{Waldron1982} for forward kinematics
calculation of serial manipulators, and in \cite{AngelesLee1988,Angeles2003}
to compute the motion equations respectively the inverse dynamics solution.
It is used in many recursive $O\left( n\right) $ forward dynamics algorithms
such as \cite%
{BaeHwangHaug1988,LillyOrin1991,Naudet2003,Orin1979,Rodriguez1991} where the
relations (\ref{VhRec}) and (\ref{Jh}) play a central role. In the so-called
'spatial operator algebra' \cite{Rodriguez1991}, hybrid screw entities are
called 'spatial vectors'. The hybrid form is deemed computationally
efficient since the transformations only involve translations. The actual
configuration of the chain is not discussed in these publications, but it
enters via the vectors $\mathbf{e}_{i}\left( \mathbf{q}\right) $ and $%
\mathbf{r}_{i}\left( \mathbf{q}\right) $, respectively $\mathbf{d}%
_{i,j}\left( \mathbf{q}\right) $. In \cite{LillyOrin1991} the inverse
transformation $\mathbf{Ad}_{\mathbf{r}_{i,j}}^{-1}$ was denoted with $^{j}%
\mathbf{X}_{i}$ (not to be confused with (\ref{Xb})), and the screw
coordinate vector ${^{0}}\mathbf{X}_{j}^{j}$ in (\ref{Xh}) with $\phi _{j}$.
In \cite{Rodriguez1991}, $\mathbf{Ad}_{\mathbf{r}_{i,j}}^{-1}$ was denoted
with $\phi _{i,j}^{T}$, and ${^{0}}\mathbf{X}_{j}^{j}$ with $\mathbf{H}%
_{j}^{T}$. The transposed matrices appear since they arose from the
transformation of wrenches.
\end{remark}

\subsubsection{Hybrid System Jacobian and its Decomposition}

The \emph{hybrid system Jacobian}, which determines the overall hybrid twist
vector $\mathsf{V}^{\text{h}}=(\mathbf{V}_{1}^{\text{h}},\ldots ,\mathbf{V}%
_{n}^{\text{h}})^{T}$ according to%
\begin{equation}
\mathsf{V}^{\text{h}}=\mathsf{J}^{\text{h}}\dot{\mathbf{q}}=\mathsf{A}^{%
\text{h}}\mathsf{X}^{\text{h}}\dot{\mathbf{q}},  \label{Vsb2}
\end{equation}%
is decomposed in terms of%
\begin{equation}
\mathsf{A}^{\text{h}}:=\left( 
\begin{array}{ccccc}
\mathbf{I} & \mathbf{0} & \mathbf{0} &  & \mathbf{0} \\ 
\mathbf{Ad}_{\mathbf{r}_{2,1}} & \mathbf{I} & \mathbf{0} &  & \mathbf{0} \\ 
\mathbf{Ad}_{\mathbf{r}_{3,1}} & \mathbf{Ad}_{\mathbf{r}_{3,2}} & \mathbf{I}
&  & \mathbf{0} \\ 
\vdots & \vdots & \ddots & \ddots &  \\ 
\mathbf{Ad}_{\mathbf{r}_{n,1}} & \mathbf{Ad}_{\mathbf{r}_{n,2}} & \cdots & 
\mathbf{Ad}_{\mathbf{r}_{n,n-1}} & \mathbf{I}%
\end{array}%
\right) ,~~~\mathsf{X}^{\text{h}}:=\mathrm{diag}~\left( {^{0}}\mathbf{X}%
_{1}^{1},\ldots ,{^{0}}\mathbf{X}_{n}^{n}\right) .  \label{Asb}
\end{equation}%
In analogy to (\ref{Invsb}), $\mathsf{A}^{\text{h}}$ can be resolved as
power series using the relation $\mathsf{A}^{\text{h}}=(\mathbf{I}-\mathsf{T}%
^{\text{h}})^{-1}$ with the $6n\times 6n$ matrix%
\begin{equation}
\mathsf{T}^{\text{h}}:=\left( 
\begin{array}{ccccc}
\mathbf{0} & \mathbf{0} & \mathbf{0} &  & \mathbf{0} \\ 
\mathbf{Ad}_{\mathbf{r}_{2,1}} & \mathbf{0} & \mathbf{0} &  &  \\ 
\mathbf{0} & \mathbf{Ad}_{\mathbf{r}_{3,2}} & \mathbf{0} &  &  \\ 
\vdots & \vdots & \ddots & \ddots &  \\ 
\mathbf{0} & \mathbf{0} & \cdots & \mathbf{Ad}_{\mathbf{r}_{n,n-1}} & 
\mathbf{0}%
\end{array}%
\right) .
\end{equation}%
This leads to the inverse $(\mathsf{A}^{\text{h}})^{-1}=(\mathbf{I}-\mathsf{T%
}^{\text{h}})$, and a solution $\dot{\mathbf{q}}$ of (\ref{Vsb2}) of the
form (\ref{solqdot1}).

\subsection{Mixed Form of Rigid Body Twists}

\subsubsection{Mixed Twists}

When formulating the Newton-Euler equations of rigid bodies, it can be
beneficial to use the body-fixed angular velocity and the translational
velocity measured at the body-fixed BFR but resolved in the IFR. This is
called the \emph{mixed twist} denoted with $\mathbf{V}_{i}^{\text{m}}=(%
\mathbf{\omega }_{i}^{\text{b}},\dot{\mathbf{r}}_{i})^{T}$. It is used in
MBS dynamics modeling \cite{Shabana}, basically because when using the mixed
twist the Newton-Euler equations w.r.t. the COM are decoupled, and because
the body-fixed inertia tensor is constant (see also the companion paper \cite%
{Part2}). The mixed twist is readily found as%
\begin{equation}
\mathbf{V}_{i}^{\text{m}}=\dot{q}_{1}\left( 
\begin{array}{c}
{^{i}}\mathbf{e}_{1} \\ 
\mathbf{b}_{i,1}\times \mathbf{e}_{1}+\mathbf{e}_{1}h_{1}%
\end{array}%
\right) +\dot{q}_{2}\left( 
\begin{array}{c}
{^{i}}\mathbf{e}_{2} \\ 
\mathbf{b}_{i,2}\times \mathbf{e}_{2}+\mathbf{e}_{2}h_{2}%
\end{array}%
\right) +\ldots +\dot{q}_{i}\left( 
\begin{array}{c}
{^{i}}\mathbf{e}_{i} \\ 
\mathbf{b}_{i,i}\times \mathbf{e}_{i}+\mathbf{e}_{i}h_{i}%
\end{array}%
\right) .  \label{Vmi}
\end{equation}%
As in (\ref{Vbsi}), $\mathbf{e}_{j}$ is a unit vector along the axis of
joint $j$ measured and resolved in the IFR $\mathcal{F}_{0}$, and $\mathbf{b}%
_{i,j}$ is the position vector of a point on the axis measured in the BFR $%
\mathcal{F}_{i}$ of body $i$ and resolved in the IFR. The mixed twist is
related to the body-fixed, spatial, and hybrid form via%
\begin{equation}
\mathbf{V}_{i}^{\text{m}}=\left( 
\begin{array}{cc}
\mathbf{I} & \mathbf{0} \\ 
\mathbf{0} & \mathbf{R}_{i}%
\end{array}%
\right) \mathbf{V}_{i}^{\text{b}}=\left( 
\begin{array}{cc}
\mathbf{R}_{i}^{T} & \mathbf{0} \\ 
\mathbf{0} & \mathbf{I}%
\end{array}%
\right) \mathbf{V}_{i}^{\text{h}}=\left( 
\begin{array}{cc}
\mathbf{R}_{i}^{T} & \mathbf{0} \\ 
-\widetilde{\mathbf{r}}_{i} & \mathbf{I}%
\end{array}%
\right) \mathbf{V}_{i}^{\text{s}}.  \label{Vmbhs}
\end{equation}

\subsubsection{Mixed Jacobian and Recursive Relations}

The expression (\ref{Vmi}) is written as%
\begin{equation}
\mathbf{V}_{i}^{\text{m}}=\sum_{j\leq i}\mathbf{J}_{i,j}^{\text{m}}\dot{q}%
_{j}=\mathsf{J}_{i}^{\text{m}}\dot{\mathbf{q}}
\end{equation}%
where the \emph{mixed Jacobian} of body $i$ is introduced as%
\begin{equation}
\mathsf{J}_{i}^{\text{m}}=\left( \mathbf{J}_{i,1}^{\text{m}},\ldots ,\mathbf{%
J}_{i,i}^{\text{m}},\mathbf{0}\ldots ,\mathbf{0}\right) .
\end{equation}%
The elements in the instantaneous joint screw coordinate vectors $\mathbf{J}%
_{i,j}^{\text{m}}$ in (\ref{Vmi}) are not consistently resolved in one
frame. Rather ${^{i}}\mathbf{e}_{j}$ is resolved in BFR $\mathcal{F}_{i}$
and $\mathbf{e}_{j}$ in the IFR. The mixed Jacobian can thus not be derived
via frame transformations. Starting from the body-fixed Jacobian, yields%
\begin{equation}
\begin{tabular}{|llll|}
\hline
&  &  &  \\ 
& $\mathbf{J}_{i,j}^{\text{m}}%
\hspace{-2ex}%
$ & $=\left( 
\begin{array}{cc}
\mathbf{R}_{i}^{T} & \mathbf{0} \\ 
\widetilde{\mathbf{r}}_{i,j} & \mathbf{I}%
\end{array}%
\right) {^{0}}\mathbf{X}_{j}^{j},\ \ j\leq i$ &  \\ 
&  &  &  \\ \hline
\end{tabular}
\label{Jm}
\end{equation}%
where ${^{0}}\mathbf{X}_{j}^{j}$ are the screw coordinates of joint $j$
measured in frame $\mathcal{F}_{j}$ and resolved in the IFR, given in (\ref%
{Xh}). The difference to (\ref{Jh}) is that the angular and translational
part are resolved in different frames. The expression (\ref{Jm}) can be
written in recursive form%
\begin{equation}
\mathbf{J}_{i,j}^{\text{m}}=\left\{ 
\begin{array}{cl}
\left( 
\begin{array}{cc}
\mathbf{R}_{i,i-1} & \mathbf{0} \\ 
\widetilde{\mathbf{r}}_{i,i-1}\mathbf{R}_{i} & \mathbf{I}%
\end{array}%
\right) \mathbf{J}_{i-1,j}^{\text{h}}, & j<i%
\vspace{-1ex}
\\ 
&  \\ 
\left( 
\begin{array}{cc}
\mathbf{R}_{i}^{T} & \mathbf{0} \\ 
\mathbf{0} & \mathbf{I}%
\end{array}%
\right) {^{0}}\mathbf{X}_{i}^{i}, & j=i\ .%
\end{array}%
\right. \ \   \label{Jmrec}
\end{equation}

This directly translates to a recursive relation for the mixed twists within
a kinematic chain%
\begin{equation}
\mathbf{V}_{i}^{\text{m}}=\left( 
\begin{array}{cc}
\mathbf{R}_{i,i-1} & \mathbf{0} \\ 
\widetilde{\mathbf{r}}_{i,i-1}\mathbf{R}_{i} & \mathbf{I}%
\end{array}%
\right) \mathbf{V}_{i-1}^{\text{m}}+\mathbf{J}{_{i,i}^{\text{m}}}\dot{q}_{i}.
\label{Vmrec}
\end{equation}

\subsubsection{Mixed System Jacobian and its Decomposition}

The overall mixed twist vector $\mathsf{V}^{\text{m}}=(\mathbf{V}_{1}^{\text{%
m}},\ldots ,\mathbf{V}_{n}^{\text{m}})^{T}$ can be expressed in terms of the
system Jacobian $\mathsf{J}^{\text{m}}$ as%
\begin{equation}
\mathsf{V}^{\text{m}}=\mathsf{J}^{\text{m}}\dot{\mathbf{q}}=\mathsf{A}^{%
\text{m}}\mathsf{X}^{\text{m}}\dot{\mathbf{q}}
\end{equation}%
with $\mathsf{X}^{\text{m}}:=\mathsf{X}^{\text{h}}$ and the matrix $\mathsf{A%
}^{\text{m}}$ as in (\ref{Asb}) but with the $\mathbf{Ad}_{\mathbf{r}_{i,j}}$
replaced by the matrix in (\ref{Jm}). This allows for a closed form
inversion of $\mathsf{A}^{\text{m}}$ analogous to that of $\mathsf{A}^{\text{%
h}}$.

\subsection{Relation of the different Forms}

The introduced twists are related by certain (not necessarily frame)
transformations, and it is occasionally desirable to switch between them.
From their definitions in (\ref{Vbs}) it is clear that body-fixed and
spatial twists, and thus the corresponding Jacobians, are related by%
\begin{equation}
\mathbf{V}_{i}^{\text{s}}=\mathbf{Ad}_{\mathbf{C}_{i}}\mathbf{V}_{i}^{\text{b%
}},\ \ \ \ \ \mathsf{J}_{i}^{\text{s}}=\mathbf{Ad}_{\mathbf{C}_{i}}\mathsf{J}%
_{i}^{\text{b}}.  \label{VsVb}
\end{equation}%
Evaluating this in the reference configuration $\mathbf{q}=\mathbf{0}$ leads
to the relation of joint screw coordinates (\ref{XY}). The body-fixed twist
is related to its hybrid version by a coordinate transformation determined
by the rotation $\mathbf{R}_{i}$ matrix, aligning the body frame with the
IFR,%
\begin{equation}
\mathbf{V}_{i}^{\text{h}}=\left( 
\begin{array}{cc}
\mathbf{R}_{i} & \mathbf{0} \\ 
\mathbf{0} & \mathbf{R}_{i}%
\end{array}%
\right) \mathbf{V}_{i}^{\text{b}}=\mathbf{Ad}_{\mathbf{R}_{i}}\mathbf{V}%
_{i}^{\text{b}}.  \label{VhVb}
\end{equation}%
The transformation (\ref{VhVb}) applies to a general hybrid screw, and in
particular to the joint screws (\ref{Xb}), (\ref{Xh}) and Jacobians (\ref%
{JbX}), (\ref{Jh}):%
\begin{equation}
{^{0}}\mathbf{X}_{j}^{j}=\mathbf{Ad}_{\mathbf{R}_{j}}{^{j}}\mathbf{X}_{j},\
\ \ \ \ \mathsf{J}_{i}^{\text{h}}=\mathbf{Ad}_{\mathbf{R}_{i}}\mathsf{J}%
_{i}^{\text{b}}.  \label{JhJb}
\end{equation}%
Combining (\ref{JhJb}) and (\ref{XY}) yields the relation of hybrid and
spatial versions of joint screws%
\begin{equation}
\mathbf{Y}_{j}=\mathbf{Ad}_{\mathbf{r}_{i}}{^{0}}\mathbf{X}_{j}^{j}
\label{YsXh}
\end{equation}%
with the current position vector $\mathbf{r}_{i}$ of body $i$ in $\mathbf{C}%
_{i}$. From (\ref{VsVb}) and (\ref{VhVb}) follows%
\begin{equation}
\mathbf{V}_{i}^{\text{s}}=\mathbf{Ad}_{\mathbf{C}_{i}}\mathbf{Ad}_{\mathbf{R}%
_{i}}^{-1}\mathbf{V}_{i}^{\text{h}}=\mathbf{Ad}_{\mathbf{r}_{i}}\mathbf{V}%
_{i}^{\text{h}}  \label{VsVh}
\end{equation}%
and thus%
\begin{equation}
\mathbf{J}_{j}^{\text{s}}=\mathbf{Ad}_{\mathbf{r}_{i}}\mathbf{J}_{i,j}^{%
\text{h}},\ \ j\leq i.  \label{JhJs}
\end{equation}%
This describes the change of reference point from the BFR of body $i$ to the
IFR. The transformations between the different forms of twists and joint
screws are summarized in table \ref{tabTransTwists}. 
\begin{table}[h] \centering%
\begin{tabular}{l|cccc|}
\cline{2-2}\cline{2-5}
& $\mathbf{V}_{i}^{\text{s}}$ & $\mathbf{V}_{i}^{\text{b}}$ & $\mathbf{V}%
_{i}^{\text{h}}$ & $\mathbf{V}_{i}^{\text{m}}$ \\ \hline
\multicolumn{1}{|l|}{$\mathbf{V}_{i}^{\text{s}}$} & $\mathbf{I}$ & $\mathbf{%
Ad}_{\mathbf{C}_{i}}$~(\ref{VsVb}) & $\mathbf{Ad}_{\mathbf{r}_{i}}$ (\ref%
{VsVh}) & (\ref{Vmbhs}) \\ 
\multicolumn{1}{|l|}{$\mathbf{V}_{i}^{\text{b}}$} & $\mathbf{Ad}_{\mathbf{C}%
_{i}}^{-1}$ & $\mathbf{I}$ & $\mathbf{Ad}_{\mathbf{R}_{i}^{T}}$ & (\ref%
{Vmbhs}) \\ 
\multicolumn{1}{|l|}{$\mathbf{V}_{i}^{\text{h}}$} & $\mathbf{Ad}_{-\mathbf{r}%
_{i}}$ & $\mathbf{Ad}_{\mathbf{R}_{i}}$ (\ref{VhVb}) & $\mathbf{I}$ & (\ref%
{Vmbhs}) \\ 
\multicolumn{1}{|l|}{$\mathbf{V}_{i}^{\text{m}}$} & (\ref{Vmbhs}) & (\ref%
{Vmbhs}) & (\ref{Vmbhs}) & $\mathbf{I}$ \\ \hline
\end{tabular}%
\ \ \ \ \ \ \ \ \ \ \ \ \ \ \ \ \ 
\begin{tabular}{c|ccc|}
\cline{2-4}
& $\mathbf{Y}_{i}$ & ${^{0}}\mathbf{X}_{i}^{i}$ & ${^{i}}\mathbf{X}_{i}$ \\ 
\hline
\multicolumn{1}{|c|}{$\mathbf{Y}_{i}$} & $\mathbf{I}$ & $\mathbf{Ad}_{%
\mathbf{r}_{i}}$ (\ref{YsXh}) & $\mathbf{Ad}_{\mathbf{A}_{i}}$ (\ref{XY}) \\ 
\multicolumn{1}{|c|}{${^{0}}\mathbf{X}_{i}^{i}$} & $\mathbf{Ad}_{-\mathbf{r}%
_{i}}$ & $\mathbf{I}$ & $\mathbf{Ad}_{\mathbf{R}_{i}}$ (\ref{JhJb}) \\ 
\multicolumn{1}{|c|}{${^{i}}\mathbf{X}_{i}$} & $\mathbf{Ad}_{\mathbf{A}%
_{i}}^{-1}$ & $\mathbf{Ad}_{\mathbf{R}_{i}^{T}}$ & $\mathbf{I}$ \\ \hline
\end{tabular}%
\caption{Transformation of the different representations of twists and joint
screw coordinates}\label{tabTransTwists}%
\end{table}%

It should be finally mentioned that the screw coordinates ${^{i}}\mathbf{X}%
_{i}$ and ${^{0}}\mathbf{X}_{i}^{i}$ are just different coordinates for the
same geometric object, namely of the instantaneous joint screw of joint $i$
measured in the BFR at body $i$ but either resolved in this BFR or in the
IFR. The vector $\mathbf{Y}_{i}$ on the other hand are a snapshot of the
joint screw coordinates of joint $i$ in spatial representation at the
reference $\mathbf{q}=\mathbf{0}$.

\begin{remark}[Computational effiency]
It is clear from (\ref{VbRec}), (\ref{Vsrec}), (\ref{VhRec}), and (\ref%
{Vmrec}) that the number of numerical operations differ between the four
different representations of twists. This allows for selecting the most
efficient one when a kinematic analysis is envisaged. In \cite%
{OrinSchrader1984} the problem of determining the twists of the terminal
body in a kinematic chain (robot end-effector) was analyzed for body-fixed,
spatial, and hybrid form. This study suggests that the spatial
representation is computationally most efficient. A conclusive analysis of
all four forms has not yet been reported. Moreover, the general situation
includes the dynamic analysis. This was partly addressed in \cite%
{Stelzle1995,Yamane2009}.
\end{remark}

\subsection{Example (continued)}

The Jacobian in body-fixed and spatial representation is determined for the
example in section \ref{secExamp1}. The instantaneous screw coordinates in
body-fixed representation are readily found with (\ref{JbX}). For instance,
the instantaneous screw coordinate vector of joint 1 expressed in the
body-fixed frame on body 3 is%
\begin{eqnarray*}
\mathbf{J}_{3,1}^{\text{b}} &=&\mathbf{Ad}_{\mathbf{C}_{3,1}\mathbf{A}%
_{1}^{-1}}\mathbf{Y}_{1}=\mathbf{Ad}_{\mathbf{C}_{3,1}}{^{1}}\mathbf{X}_{1}
\\
&=&\left(
0,0,1,(d_{2}+d_{3})s_{2}+(x_{1}+x_{3}-d_{3})s_{23},d_{2}-x_{1}-(d_{2}+d_{3})c_{2}+(d_{3}-x_{1}-x_{3})c_{23},0\right) ^{T}.
\end{eqnarray*}%
Proceeding analogously for the other joints, yields the body-fixed Jacobian
of body 3 as%
\begin{equation*}
\mathsf{J}_{3}^{\text{b}}\left( \mathbf{q}\right) =\left( 
\begin{array}{ccccc}
0 & 0 & 0 & 0 & 0 \\ 
0 & 0 & 0 & 0 & 0 \\ 
1 & 1 & 1 & 0 & 0 \\ 
({d_{2}}+{d_{3}})s_{2}+({x}_{{1}}+{x_{3}}-{d_{3}})s_{23} & s_{3}({x_{2}}+{%
x_{3}}-{d_{2}}-{d_{3}}) & 0 & 0 & 0 \\ 
{d_{2}}-{x}_{{1}}-({d_{2}}+{d_{3}})c_{2}+({d_{3}}-{x}_{{1}}-{x_{3}})c_{23} & 
c_{3}({d_{2}}+{d_{3}}-{x_{2}}-{x_{3}})-{d_{3}}-{x_{2}} & {x_{3}}-{d_{3}} & 0
& 0 \\ 
0 & 0 & 0 & 0 & 0%
\end{array}%
\right) .
\end{equation*}%
The body-fixed twist of body 3 is therewith $\mathbf{V}_{3}^{\text{b}}=%
\mathsf{J}_{3}^{\text{b}}\left( \mathbf{q}\right) \dot{\mathbf{q}}$. Again
details for body 4 and 5 are omitted due space limitation.

The spatial representation of the screw coordinates of joint $1,\ldots ,4$,
for instance, are found with (\ref{JsX})%
\begin{eqnarray*}
\mathbf{J}_{1}^{\text{s}}\left( \mathbf{q}\right) &=&\left(
0,0,1,0,0,0\right) ^{T} \\
\mathbf{J}_{2}^{\text{s}}\left( \mathbf{q}\right) &=&\left(
0,0,1,-d_{2}s_{1},d_{2}c_{1},0\right) ^{T} \\
\mathbf{J}_{3}^{\text{s}}\left( \mathbf{q}\right) &=&\left(
0,0,1,(d_{2}+d_{3})s_{23}-d_{2}s_{1},d_{2}c_{1}-(d_{2}+d_{3})c_{12},0\right)
^{T} \\
\mathbf{J}_{4}^{\text{s}}\left( \mathbf{q}\right) &=&\frac{1}{\sqrt{2}}%
(-c_{123},-s_{123},1,(d_{2}+d_{3})s_{12}-d_{2}s_{1}+(d_{4}+h_{4}-d_{3})s_{123},
\\
&&\ \ \ \ \ \
d_{2}c_{1}-(d_{2}+d_{3})c_{12}+(d_{3}-d_{4}-h_{4})c_{123},d_{2}s_{23}-(d_{2}+d_{3})s_{3})^{T}.
\end{eqnarray*}%
That for joint 5 is omitted again. These $\mathbf{J}_{i}^{\text{s}}$
constitute the spatial Jacobians $\mathsf{J}_{i}^{\text{s}}$ in (\ref{Js}).

\section{Conclusions and Outlook}

Screw and Lie group theory gives rise to compact formulations of the
equations governing the MBS kinematics in terms of relative (joint)
coordinates. This is beneficial for the actual modeling process as well as
for the implementation of MBS algorithms and their computational properties.
The frame invariance of these concepts allows for expressing the relevant
modeling objects as suited best for a particular application. In particular
the MBS kinematics can be formulated without introduction of body-fixed
joint frames. This is a central result that gives rise to maximal
flexibility as opposed to the use of modeling conventions like
Denavit-Hartenberg parameters. These results have been published over the
last two decades, but they have not been presented within a uniform MBS
framework. In this paper, screw and Lie group theory have been employed to
provide such a framework. Decisive for the computational efficiency is the
actual representation of rigid body twists and accelerations. Four commonly
used forms were recalled, and the recursive algorithms for MBS kinematics
where presented. The corresponding recursive algorithms for evaluation of
the motion equations are presented in the accompanying paper \cite{Part2}.

The reader used to work with the classical body-fixed twists should be able
to directly apply the presented modeling paradigm for MBS kinematics using
the relation (\ref{POEY}) to determine the body configurations and (\ref{JbX}%
) to determine the Jacobian while having the freedom to choose arbitrary BFR
and IFR. This applies likewise to the spatial, hybrid, and mixed twists.

The full potential of Lie group formulations is yet to be explored in future
research. This regards the modeling steps as well as the computational
properties, in particular given a current trend in computational MBS
dynamics to put more emphasize on user friendly modeling and on tailored
simulation codes. A forthcoming paper will address MBS with general
topology. To this end, the loop closure constraints are formulated in the
form of a POE. Redundant loop constraints are still a major challenge. It is
already known that the loop constraints can be concisely formulated in terms
of joint screws, but even more that they can be reduced to a non-redundant
constraint system by means of simple operations on the joint screw system 
\cite{ConstraintReduct2014}.

\section*{A Rigid Body Motions and the Lie Group $SE\left( 3\right) $}

For an introduction to screws and to the motion Lie group $SE\left( 3\right) 
$ the reader is referred to the text books \cite%
{Angeles2003,LynchPark2017,Murray,Selig}.

\subsection*{A.1 Finite Rigid Body Motions as Frame Transformations -- $%
SE\left( 3\right) $}

A frame $\{O_{i};\vec{e}_{i,1},\vec{e}_{i,2},\vec{e}_{i,3}\}$ consists of a
point $O_{i}\in E^{3}$ (its origin) and a basis triad $\{\vec{e}_{i,1},\vec{e%
}_{i,2},\vec{e}_{i,3}\}$, with $\vec{e}_{i,k}\in {\mathbb{E}}^{3}$, in which
vectors are resolved. A \emph{change of basis} from $\{\vec{e}_{i,1},\vec{e}%
_{i,2},\vec{e}_{i,3}\}$ to $\{\vec{e}_{j,1},\vec{e}_{j,2},\vec{e}_{j,3}\}$
is a is a coordinate transformation from coordinates resolved in $\{\vec{e}%
_{i,1},\vec{e}_{i,2},\vec{e}_{i,3}\}$ to coordinates resolved in $\{\vec{e}%
_{j,1},\vec{e}_{j,2},\vec{e}_{j,3}\}$. A \emph{frame transformation} from $%
\{O_{i};\vec{e}_{i,1},\vec{e}_{i,2},\vec{e}_{i,3}\}$ to $\{O_{j};\vec{e}%
_{j,1},\vec{e}_{j,2},\vec{e}_{j,3}\}$ is a coordinate transformation
together with a change of origin from $O_{i}$ to $O_{j}$. When a vector is
resolved in a frame according to $\vec{r}={^{i}}r_{1}\vec{e}_{i,1}+{^{i}}%
r_{2}\vec{e}_{i,2}+{^{i}}r_{3}\vec{e}_{i,3}\}$ its component vector is
doneted by ${^{i}}\mathbf{r}=\left( {^{i}}r_{1},{^{i}}r_{2},{^{i}}%
r_{3}\right) ^{T}\in {\mathbb{R}}^{3}$. The leading superscript indicates
the frame in which it is resolved.

A rigid body is kinematically represented by a body-fixed reference frame
(BFR). Denote the BFR of body $i$ with $\mathcal{F}_{i}=\{\Omega _{i},\vec{e}%
_{i,1},\vec{e}_{i,2},\vec{e}_{i,3}\}$. Its motion is thus described as the
relative motion of the BFR w.r.t. a global inertial frame (IFR) $\mathcal{F}%
_{0}=\{O,\vec{e}_{1},\vec{e}_{2},\vec{e}_{3}\}$. The location of $\mathcal{F}%
_{i}$ is described by its global position vector $\vec{r}=\overline{O\Omega
_{i}}$. When resolved in the IFR $\mathcal{F}_{0}$, its coordinate vector is 
${^{0}}\mathbf{r}_{i}\in \mathbb{R}^{3}$. The orientation is described by a
rotation matrix $\mathbf{R}_{0,i}\in SO\left( 3\right) $ that transforms
coordinates of a vector $\mathbf{x}$ resolved in the BFR to its coordinates
when resolved in the IFR according to ${^{0}}\mathbf{x}=\mathbf{R}_{0,i}{^{i}%
}\mathbf{x}$. In the following the index 0 is omitted, i.e. $\mathbf{x}$ is
the coordinate vector resolved in the IFR.

If ${^{i}}\mathbf{b}\in \mathbb{R}^{3}$ is the position vector of a point $P$
of the body resolved in the BFR, the position vector of point $P$ measured
and resolved in IFR is $\mathbf{s}=\mathbf{r}+\mathbf{R}_{i}{^{i}}\mathbf{b}$%
. This transformation can be written compactly using homogenous point
coordinates%
\begin{equation}
\left( 
\begin{array}{c}
\mathbf{s} \\ 
\ 1%
\end{array}%
\right) =\left( 
\begin{array}{cc}
\mathbf{R}_{i} & \mathbf{r}_{i} \\ 
\mathbf{0} & 1%
\end{array}%
\right) \left( 
\begin{array}{c}
{^{i}}\mathbf{b} \\ 
\ 1%
\end{array}%
\right) .  \label{homoTrans}
\end{equation}%
This is a frame transformation, i.e. it describes the transformation due to
the rotation as well as due to the displacement of the origin of the
reference frame. As this holds for any point of the rigid body, the matrix%
\begin{equation}
\mathbf{C}_{i}=\left( 
\begin{array}{cc}
\mathbf{R}_{i} & \mathbf{r}_{i} \\ 
\mathbf{0} & 1%
\end{array}%
\right) \in SE\left( 3\right)  \label{C}
\end{equation}%
describes the configuration of the BFR $\mathcal{F}_{i}$ w.r.t. to the IFR,
which is referred to as the \emph{absolute configuration} of $\mathcal{F}%
_{i} $ (as it refers to the global IFR). For simplicity the configuration is
alternatively denoted by the pair $C_{i}=\left( \mathbf{R}_{i},\mathbf{r}%
_{i}\right) $. $SE\left( 3\right) $ is the group of isometric orientation
preserving transformations of Euclidean spaces. It is is commonly
represented as matrix group with elements as in (\ref{homoTrans}). The
inverse of the transformation (\ref{C}) is 
\begin{equation}
\mathbf{C}_{i}^{-1}=\left( 
\begin{array}{cc}
\mathbf{R}_{i}^{T} & -\mathbf{R}_{i}^{T}\mathbf{r}_{i} \\ 
\mathbf{0} & 1%
\end{array}%
\right) =\left( 
\begin{array}{cc}
\mathbf{R}_{i} & -\mathbf{R}_{i}\mathbf{r}_{i} \\ 
\mathbf{0} & 1%
\end{array}%
\right)
\end{equation}%
respectively $C_{i}^{-1}=\left( \mathbf{R}_{i}^{T},-\mathbf{R}_{i}^{T}%
\mathbf{r}_{i}\right) $. Let $C^{\prime }$ and $C^{\prime \prime }$ be two
frame transformations. The product $C^{\prime }\cdot C^{\prime \prime
}=\left( \mathbf{R}^{\prime }\mathbf{R}^{\prime \prime },\mathbf{r}^{\prime
}+\mathbf{R}^{\prime }\mathbf{r}^{\prime \prime }\right) $, respectively $%
\mathbf{C}^{\prime }\cdot \mathbf{C}^{\prime \prime }$, describes the
overall frame transformation.

Now consider two bodies, i.e. two RFRs $\mathcal{F}_{i}$ and $\mathcal{F}%
_{j} $, and denote their absolute configuration with $\mathbf{C}_{i}$ and $%
\mathbf{C}_{j}$, respectively. The \emph{relative configuration} of body $j$
w.r.t. body $i$ is 
\begin{equation}
\mathbf{C}_{i,j}=\mathbf{C}_{i}^{-1}\mathbf{C}_{j}=\left( 
\begin{array}{cc}
\mathbf{R}_{i}^{T}\mathbf{R}_{j} & \ \ \mathbf{R}_{i}^{T}\left( \mathbf{r}%
_{j}-\mathbf{r}_{i}\right) \\ 
\mathbf{0} & 1%
\end{array}%
\right) =\left( 
\begin{array}{cc}
\mathbf{R}_{i,j} & {^{i}}\mathbf{r}_{i,j} \\ 
\mathbf{0} & 1%
\end{array}%
\right)  \label{Cij}
\end{equation}%
with the relative rotation matrix $\mathbf{R}_{i,j}$ and the relative
displacement vector $\mathbf{r}_{i,j}:=\mathbf{r}_{j}-\mathbf{r}_{i}$
resolved in the RFR $\mathcal{F}_{i}$ on body $i$. The configuration of body 
$j$ is then expressed in terms of the configuration of body $i$ and the
relative configuration as $\mathbf{C}_{j}=\mathbf{C}_{i}\mathbf{C}_{i,j}$.
Analogously, $\mathbf{C}_{j,i}=\mathbf{C}_{j}^{-1}\mathbf{C}_{i}$ is the
relative configuration of body $i$ w.r.t. body $j$. Clearly, $\mathbf{C}%
_{j,i}=\mathbf{C}_{i,j}^{-1}$. As special case, the absolute configuration
of body $i$ is $\mathbf{C}_{i}=\mathbf{C}_{0,i}=\mathbf{C}_{0}^{-1}\mathbf{C}%
_{i}$ with $\mathbf{C}_{0}=\mathbf{I}$.

Throughout the paper, an $SE\left( 3\right) $ matrix is always considered to
represent the relative configuration of two frames (also $\mathbf{C}_{i}$ is
the configuration of body $i$ relative to the IFR).

Rotation matrices form the 3-dimensional special orthogonal group, denoted $%
SO\left( 3\right) $. This is a Lie group, which means that for any rotation
matrix there is a unique inverse and that there is a smooth
parameterization. A smooth canonical parameterization is the description in
terms of axis and angle. Moreover, Euler's theorem states that any finite
rotation can be achieved by a rotation about an axis. Let $\mathbf{e}\in {%
\mathbb{R}}^{3}$ be unit vector along the rotation axis, $\varphi $ the
rotation angle, and denote with $\mathbf{\xi }:=\mathbf{e}\varphi $ is the
'scaled rotation axis'. The corresponding rotation matrix is $\mathbf{R}%
\left( \varphi ,\mathbf{e}\right) =\exp \widetilde{\mathbf{\xi }}$ with%
\begin{eqnarray}
\exp \widetilde{\mathbf{\xi }} &=&\mathbf{I}+\tfrac{\sin \left\Vert \mathbf{%
\xi }\right\Vert \,}{\left\Vert \mathbf{\xi }\right\Vert }\widetilde{\mathbf{%
\xi }}+\tfrac{1-\cos \left\Vert \mathbf{\xi }\right\Vert }{\left\Vert 
\mathbf{\xi }\right\Vert ^{2}}\,\widetilde{\mathbf{\xi }}^{2}  \label{SO3exp}
\\
&=&\mathbf{I}+\mathrm{sinc}\left\Vert \mathbf{\xi }\right\Vert \widetilde{%
\mathbf{\xi }}+\tfrac{1}{2}\mathrm{sinc}^{2}\tfrac{\left\Vert \mathbf{\xi }%
\right\Vert }{2}\,\widetilde{\mathbf{\xi }}^{2}  \label{SO3exp2} \\
&=&\mathbf{I}+\sin \varphi \widetilde{\mathbf{e}}+\left( 1-\cos \varphi
\right) \,\widetilde{\mathbf{e}}^{2}  \label{SO3exp3}
\end{eqnarray}%
where $\widetilde{\mathbf{\xi }}\in so\left( 3\right) $ denotes the skew
symmetric matrix associated to the vector $\mathbf{\xi }$. This is known as
the Euler-Rodrigues formula \cite{Angeles2003,McCarthy1990}. The rotation
matrix (\ref{SO3exp}) takes a frame from its initial to its final
orientation, where the rotation vector $\mathbf{e}$ is resolved in its
initial orientation. If the rotation matrix is considered to transform
coordinates from different frames, $\mathcal{F}_{i}$ and $\mathcal{F}_{j}$,
then $\mathbf{R}_{i,j}=\exp {^{i}}\widetilde{\mathbf{\xi }}$ where the axis
is resolved in $\mathcal{F}_{i}$.

The important concept here is that of the exp mapping. In (\ref{SO3exp})
this is the matrix exponential mapping a skew symmetric matrix to a rotation
matrix. Moreover, skew symmetric $3\times 3$ matrices form the Lie algebra
denoted $so\left( 3\right) $. Being a Lie algebra, it is a vector space
equipped with the Lie bracket, in this case the matrix commutator $\left[ 
\mathbf{\widetilde{\mathbf{x}}},\mathbf{\widetilde{\mathbf{y}}}\right] =%
\mathbf{\widetilde{\mathbf{x}}\widetilde{\mathbf{y}}-\widetilde{\mathbf{y}}%
\widetilde{\mathbf{x}}}$. This can be expressed $\left[ \mathbf{\widetilde{%
\mathbf{x}}},\mathbf{\widetilde{\mathbf{y}}}\right] =\widetilde{\mathbf{x}%
\times \mathbf{y}}$, which means that $so\left( 3\right) $ is isomorphic to
the vector space ${\mathbb{R}}^{3}$ equipped with the cross product.
Elements of $so\left( 3\right) $ can thus be represented as 3-vectors via $%
\mathbf{\widetilde{\mathbf{x}}}\in so\left( 3\right) \leftrightarrow \mathbf{%
x}\in {\mathbb{R}}^{3}.$

Frame transformations form the special Euclidean group $SE\left( 3\right)
=SO\left( 3\right) \ltimes \mathbb{R}^{3}$ --more precisely the group of
isometric and orientation preserving transformations of Euclidean spaces. A
typical element is represented as a matrix (\ref{C}). The multiplication of
these matrices, respectively the transformation (\ref{homoTrans}), reveals
that $SE\left( 3\right) $ is the semidirect product of the rotation group $%
SO\left( 3\right) $ and the translation group ${\mathbb{R}}^{3}$. This may
not seem important but it has consequences for the constraint satisfaction
when integrating MBS models described in absolute coordinates \cite%
{MMT2014,BIT}. should be $SE\left( 3\right) $ is a 6-dimensional Lie group,
so that it possesses a smooth parameterization, and for any element there is
an inverse. Chasles' theorem \cite{Angeles2003,Murray,Selig} states that any
finite rigid body displacement (i.e. frame transformation) can be achieved
by a screw motion, i.e. a rotation about a constant axis together with a
translation along this axis that is determined by the pitch.

\begin{figure}[b]
\centerline{
\includegraphics[width=9cm]{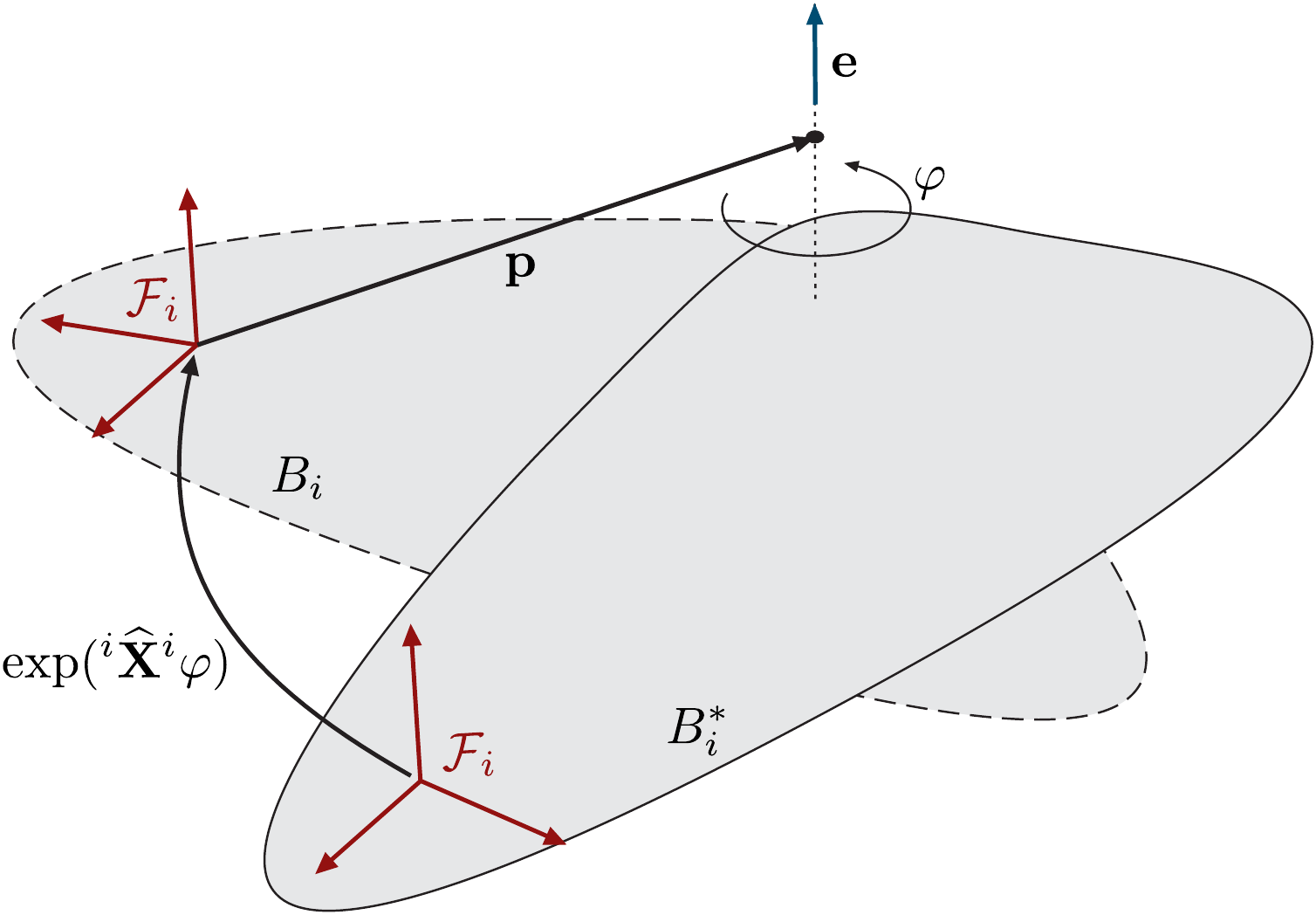}}
\caption{Screw motion of a rigid body represented by a body-fixed frame $%
\mathcal{F}_{i}$}
\label{figScrewMotion}
\end{figure}

Consider a body performing a screw motion. Let ${^{i}}\mathbf{e}$ be the
unit vector along the screw axis and $\mathbf{p}$ be the position vector of
an arbitrary point on that axis (fig. \ref{figScrewMotion}), both resolved
in the body-fixed frame $\mathcal{F}_{i}$. The vector of \emph{unit screw
coordinates} corresponding to this motion is%
\begin{equation}
{^{i}}\mathbf{X}{^{i}}=\left( 
\begin{array}{c}
{^{i}}\mathbf{e} \\ 
{^{i}}\mathbf{p}\times {^{i}}\mathbf{e}+{^{i}}\mathbf{e}h%
\end{array}%
\right)  \label{iXi}
\end{equation}%
where $h$ is the pitch of the screw \cite{BottemaRoth1990,Murray,CND2017}.
In classical screw theory literature screws are denoted by the symbol '\$'.
The vector $\left( {^{i}}\mathbf{e},{^{i}}\mathbf{p}\times {^{i}}\mathbf{e}%
\right) ^{T}$, i.e. setting $h=0$, are Pl\"{u}cker coordinates of the line
along the screw axis. Geometrically, a screw is determined by the Pl\"{u}%
cker coordinates of the line along the screw axis and the pitch. 

The screw motion, with rotation angle $\varphi $, taking the body from its
initial configuration $B_{i}$ to the final configuration $B_{i}^{\ast }$ is
determined by the matrix exponential $\exp ({^{i}}\widehat{\mathbf{X}}{^{i}}%
\varphi )$ of the matrix ${^{i}}\widehat{\mathbf{X}}{^{i}}$ defined as 
\begin{equation}
\mathbf{X}=\left( 
\begin{array}{c}
\mathbf{e} \\ 
\mathbf{p}\times \mathbf{e}+\mathbf{e}h%
\end{array}%
\right) \ \in {\mathbb{R}}^{6}\ \leftrightarrow \ \ \widehat{\mathbf{X}}%
=\left( 
\begin{array}{cc}
\widetilde{\mathbf{e}} & \ \ \ \mathbf{p}\times \mathbf{e}+\mathbf{e}h \\ 
\mathbf{0} & 0%
\end{array}%
\right) \in se\left( 3\right) .  \label{Xhat}
\end{equation}%
The exponential mapping admits the closed form expression%
\begin{equation}
\exp (\varphi \widehat{\mathbf{X}})=\left( 
\begin{array}{cc}
\exp (\varphi \widetilde{\mathbf{e}}) & \ \ \ (\mathbf{I}-\exp (\varphi 
\widetilde{\mathbf{e}}))\mathbf{p}+\varphi h\mathbf{e} \\ 
\mathbf{0} & 1%
\end{array}%
\right) ,\ \text{for }h\neq \infty  \label{SE3expP}
\end{equation}%
where $\varphi $ is the rotation angle, and for pure translations, i.e.
infinite pitch,%
\begin{equation}
\exp (\varphi \widehat{\mathbf{X}})=\left( 
\begin{array}{cc}
\mathbf{I} & \ \varphi \mathbf{e} \\ 
\mathbf{0} & 1%
\end{array}%
\right) ,\ \text{for }h=\infty
\end{equation}%
where $\varphi $ is the translation variable. The exp mapping describes the
frame transformation of a body-fixed frame from its final configuration to
its initial configuration due to a screw motion where the screw coordinates
are represented in the initial configuration. This is applicable to general
frame transformations. Let ${^{i}}\mathbf{X}{^{i}}$ be the screw coordinates
associated to the relative screw motion of body $j$ w.r.t. body $i$. The
configuration of the body-fixed frame $\mathcal{F}_{j}$ relative to $%
\mathcal{F}_{i}$ is $\mathbf{C}_{i,j}=\exp ({^{i}}\widehat{\mathbf{X}}{^{i}}%
\varphi )$, assuming that $\mathcal{F}_{j}$ and $\mathcal{F}_{i}$ initially
overlap.

Notice the direction in which the frame transformation is indicated in fig. %
\ref{figScrewMotion}. The arc points toward the frame in which the screw
coordinates are expressed.

Matrices of the form shown in (\ref{Xhat}) play obviously a key role as they
give rise to frame transformations via the exp mapping. They form the Lie
algebra $se\left( 3\right) $. To any matrix $\widehat{\mathbf{X}}\in
se\left( 3\right) $ corresponds a unique vector $\mathbf{X}=(\mathbf{\xi },%
\mathbf{\eta })^{T}\in {\mathbb{R}}^{6}$ via the isomorphism (\ref{Xhat}).
For simplicity also the notation $\exp \left( \mathbf{X}\varphi \right) $ is
used instead of $\exp (\widehat{\mathbf{X}}\varphi )$.

\begin{figure}[b]
\centerline{
\includegraphics[width=8cm]{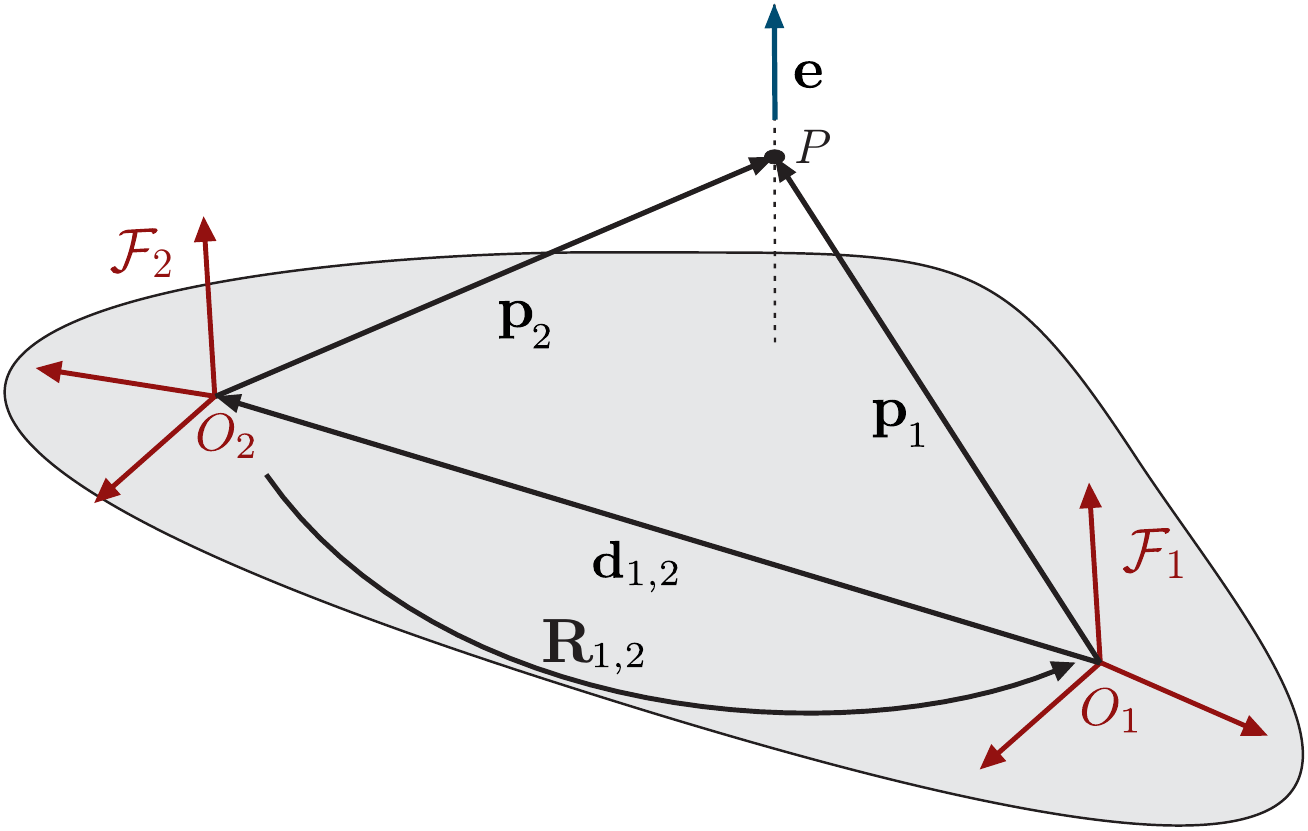}}
\caption{Frame transformation of screw coordinates according to $%
S_{1,2}=\left( \mathbf{R}_{1,2},{^{1}}\mathbf{d}_{1,2}\right) $.}
\label{figScrewTrans}
\end{figure}

Screws are geometric objects, thus frame invariant, and can be represented
in any reference frame. Consider two frames $\mathcal{F}_{1}$ and $\mathcal{F%
}_{2}$ (fig. \ref{figScrewTrans}), and let ${^{2}}\mathbf{X}{^{2}}=\left( {%
^{2}}\mathbf{e},{^{2}}\mathbf{p}_{2}\times {^{2}}\mathbf{e}+{^{2}}\mathbf{e}%
h\right) ^{T}$ be screw coordinates measured and resolved in frame $\mathcal{%
F}_{2}$. That is, ${^{2}}\mathbf{e}$ is the unit vector along the axis
resolved in $\mathcal{F}_{2}$, and ${^{2}}\mathbf{p}_{2}$ is the position
vector $\overline{O_{2}P}$ of point $P$ on the axis resolved in $\mathcal{F}%
_{2}$. Let $S_{1,2}=\left( \mathbf{R}_{1,2},{^{1}}\mathbf{d}_{1,2}\right) $
be the transformation from $\mathcal{F}_{2}$ to $\mathcal{F}_{1}$. Then the
screw coordinate vector measured and resolved in $\mathcal{F}_{1}$ is
determined by ${^{1}\widehat{\mathbf{X}}^{1}}=\mathbf{S}_{1,2}{^{2}\widehat{%
\mathbf{X}}^{2}}\mathbf{S}_{1,2}^{-1}=\mathrm{Ad}_{\mathbf{S}_{1,2}}({^{2}%
\widehat{\mathbf{X}}^{2}})$, which reads in vector notation%
\begin{equation}
{^{1}\mathbf{X}^{1}}=\left( 
\begin{array}{c}
{^{1}}\mathbf{e} \\ 
\ {^{1}}\mathbf{p}_{1}\times {^{1}}\mathbf{e}+{^{1}}\mathbf{e}h%
\end{array}%
\right) =\left( 
\begin{array}{cc}
\mathbf{R}_{1,2} & \mathbf{0} \\ 
{^{1}}\widetilde{\mathbf{d}}_{1,2}\mathbf{R}_{1,2}\ \  & \mathbf{R}_{1,2}%
\end{array}%
\right) \left( 
\begin{array}{c}
{^{2}}\mathbf{e} \\ 
\ {^{2}}\mathbf{p}_{2}\times {^{2}}\mathbf{e}+{^{2}}\mathbf{e}h%
\end{array}%
\right) =\mathbf{Ad}_{\mathbf{S}_{1,2}}{^{2}}\mathbf{X}{^{2}}  \label{Ad2}
\end{equation}%
where ${^{1}}\mathbf{p}_{1}=\mathbf{R}_{1,2}{^{2}}\mathbf{p}_{2}+{^{1}}%
\mathbf{d}_{1,2}$ is the position vector $\overline{O_{1}P}$ of the point $P$
on the screw axis measured and resolved in $\mathcal{F}_{1}$. The
transformation 
\begin{equation}
\mathrm{Ad}_{\mathbf{C}}({\widehat{\mathbf{X}}})=\mathbf{C}{\widehat{\mathbf{%
X}}}\mathbf{C}^{-1}  \label{Ad1}
\end{equation}%
with a general frame transformation $\mathbf{C}\in SE\left( 3\right) $ is
called the \emph{adjoined transformation} \cite{LynchPark2017,Murray,Selig}.
In vector notion this is%
\begin{equation}
\mathbf{Ad}_{\mathbf{C}}=\left( 
\begin{array}{cc}
\mathbf{R} & \mathbf{0} \\ 
\widetilde{\mathbf{r}}\mathbf{R}\ \  & \mathbf{R}%
\end{array}%
\right) .  \label{Ad}
\end{equation}%
The terminology stems from the fact that $\mathbf{C}\in SE\left( 3\right) $
describes a frame transformation, while $\mathbf{Ad}_{\mathbf{C}}$ describes
the corresponding transformation of screw coordinates that belong to $%
se\left( 3\right) $. It thus equally describes rigid body motions. It is
used for instance in \cite{Bauchau2011} where $\mathbf{Ad}_{\mathbf{C}}$ is
called the 'motion tensor' and in \cite{Borri2001a} where it is referred to
as 'configuration tensor' (and denoted with $\mathbf{C}$). If $\mathbf{C}%
_{i,j}$ is the configuration of body $j$ relative to body $i$, then $\mathbf{%
Ad}_{\mathbf{C}_{i,j}}$ transforms the coordinates of a screw represented in
BFR on body $j$ to its coordinate representation in the BFR on body $i$.

For the product of transformations it holds 
\begin{equation}
\mathbf{Ad}_{\mathbf{C}^{\prime }\mathbf{C}^{\prime \prime }}=\mathbf{Ad}_{%
\mathbf{C}^{\prime }}\mathbf{Ad}_{\mathbf{C}^{\prime \prime }}.
\end{equation}

For sake of compactness, with slight abuse of notation the following
notation is used%
\begin{equation}
\mathbf{Ad}_{\mathbf{R}}=\left( 
\begin{array}{cc}
\mathbf{R} & \mathbf{0} \\ 
\mathbf{0} & \mathbf{R}%
\end{array}%
\right) ,\ \text{for}\ \ C=\left( \mathbf{R},\mathbf{0}\right) \ \ \ \ \ \ \
\ \mathbf{Ad}_{\mathbf{r}}=\left( 
\begin{array}{cc}
\mathbf{I} & \mathbf{0} \\ 
\widetilde{\mathbf{r}} & \mathbf{I}%
\end{array}%
\right) ,\ \text{for}\ \ C=\left( \mathbf{I},\mathbf{r}\right) .
\label{AdRr}
\end{equation}%
This allows for splitting the frame transformation (\ref{Ad}) into the
change of reference point followed by the change of basis as $\mathbf{Ad}_{%
\mathbf{C}}=\mathbf{Ad}_{\mathbf{r}}\mathbf{Ad}_{\mathbf{R}}$. 

A useful relation is that, for any $\widehat{\mathbf{X}}\in se\left(
3\right) $ and $\mathbf{S}\in SE\left( 3\right) $, it is%
\begin{equation}
\mathbf{S}\exp (\widehat{\mathbf{X}})\mathbf{S}^{-1}=\exp (\mathbf{S}%
\widehat{\mathbf{X}}\mathbf{S}^{-1})=\exp (\mathrm{Ad}_{\mathbf{S}}(\widehat{%
\mathbf{X}})).  \label{Adexp}
\end{equation}%
For constant $\mathbf{X}$ the derivative of the exponential is 
\begin{equation}
\frac{\partial }{\partial t}\exp (t\widehat{\mathbf{X}})=\widehat{\mathbf{X}}%
\exp (t\widehat{\mathbf{X}})=\exp (t\widehat{\mathbf{X}})\widehat{\mathbf{X}}%
.  \label{derExpX}
\end{equation}

\paragraph{Notation:}

The coordinate representation of a screw requires 1) a reference point from
which the point on the screw axis is measured, and 2) a reference frame in
which the vectors are resolved. The reference point is commonly the origin
of a frame. These are indicated by the leading and trailing superscript,
respectively. In (\ref{Ad2}) they were identical, but in general, ${^{i}%
\mathbf{X}^{j}}=\left( {^{i}}\mathbf{e},{^{i}}\mathbf{p}_{j}\times {^{i}}%
\mathbf{e}+{^{i}}\mathbf{e}h\right) ^{T}$ is the screw coordinate vector
measured from the origin of $\mathcal{F}_{j}$ and resolved in $\mathcal{F}%
_{i}$. As apparent from (\ref{Ad2}) only screws measured and resolved in the
same frame ($i=j$) are related by frame transformations. For sake of
simplifying the notation, when $i=j$ the simplified notation ${^{i}\mathbf{X}%
}$ is used. The screw is then said to be \emph{represented} in $\mathcal{F}%
_{i}$.

\subsection*{A.2 Twists as instantaneous Screw Motions -- $se\left( 3\right) 
$%
\label{secTwists}%
}

The angular and translational velocity of a rigid body are summarized in the
vector $\mathbf{V}=\left( \mathbf{\omega },\mathbf{v}\right) ^{T}$ called
the rigid body twist. This is a screw, and can hence be written as matrix of
the form (\ref{Xhat})%
\begin{equation}
\widehat{\mathbf{V}}=\left( 
\begin{array}{cc}
\widetilde{\mathbf{\omega }} & \mathbf{v} \\ 
\mathbf{0} & 0%
\end{array}%
\right) .  \label{twistMat}
\end{equation}%
A general definition of twists requires specification of 1) the body of
which the twist is measured, 2) the point at which the velocity is measured,
and 3) the frame in which the velocity vectors are resolved. To this end,
the following notion is used:%
\begin{equation}
{^{k}\mathbf{V}_{i}^{j}=}\left( 
\begin{array}{c}
{^{k}}\mathbf{\omega }{_{i}} \\ 
{^{k}}\mathbf{v}{_{i}^{j}}%
\end{array}%
\right) ,\ \ \ \ \ \ 
\begin{array}{ll}
i & \text{- index of the 'object' of which the twist is measured} \\ 
j & \text{- index of the frame in which the twist is measured} \\ 
k & \text{- index of the frame in which the vectors are resolved.}%
\end{array}
\label{Vnotation}
\end{equation}%
This is the twist of body $i$ measured in $\mathcal{F}_{j}$ and resolved in $%
\mathcal{F}_{k}$. More precisely ${^{k}}\mathbf{\omega }{_{i}}$ is the
angular velocity of frame $\mathcal{F}_{i}$ measured and resolved in $%
\mathcal{F}_{k}$. Due to the translation invariance of the angular
velocities it is independent from $\mathcal{F}_{j}$. The vector ${^{k}}%
\mathbf{v}{_{i}^{j}}$ is the translational velocity of the point on body $i$
that is instantaneously traveling trough the origin of $\mathcal{F}_{j}$
resolved in $\mathcal{F}_{k}$.

When $j=k$, the simplified notation ${^{k}\mathbf{V}_{i}}$ is used. A twist
(and generally a screw) is said to be \emph{represented in frame }$\mathcal{F%
}_{k}$ if it is measured and resolved in this frame, e.g. ${^{k}\mathbf{V}%
_{i}}$ is represented in frame $k$. The \emph{body-fixed} and \emph{spatial}
representation of twists are most commonly used. The attribute 'body-fixed'
indicates that the frame in which the velocity is measured and resolved is
the body-fixed RFR, i.e. $i=j=k$. 'Spatial' is used when velocities are
measured and resolved in the IFR, i.e. $j=k=0$. To further simplify the
notation, the body-fixed and spatial twist of body $i$ is denoted by $%
\mathbf{V}_{i}^{\text{b}}=\left( \mathbf{\omega }_{i}^{\text{b}},\mathbf{v}%
_{i}^{\text{b}}\right) ^{T}:={^{i}\mathbf{V}_{i}^{i}}$ and $\mathbf{V}_{i}^{%
\text{s}}=\left( \mathbf{\omega }_{i}^{\text{s}},\mathbf{v}_{i}^{\text{s}%
}\right) ^{T}:={^{0}\mathbf{V}_{i}^{0}}$, respectively. The index 0 is
omitted throughout the paper. For a body whose motion is described by $%
\mathbf{C}_{i}\left( t\right) $ according to (\ref{C}) these are defined
analytically by%
\begin{equation}
\widehat{\mathbf{V}}_{i}^{\text{b}}=\left( 
\begin{array}{cc}
\widetilde{\mathbf{\omega }}_{i}^{\text{b}} & \mathbf{v}_{i}^{\text{b}} \\ 
\mathbf{0} & 0%
\end{array}%
\right) =\mathbf{C}_{i}^{-1}\dot{\mathbf{C}}_{i},\ \ \ \ \ \ \ \ \widehat{%
\mathbf{V}}_{i}^{\text{s}}=\left( 
\begin{array}{cc}
\widetilde{\mathbf{\omega }}_{i}^{\text{s}} & \mathbf{v}_{i}^{\text{s}} \\ 
\mathbf{0} & 0%
\end{array}%
\right) =\dot{\mathbf{C}}_{i}\mathbf{C}_{i}^{-1}.  \label{Vbs}
\end{equation}%
Therein $\widetilde{\mathbf{\omega }}_{i}^{\text{b}}=\mathbf{R}_{i}^{T}\dot{%
\mathbf{R}}_{i}$ and $\widetilde{\mathbf{\omega }}_{i}^{\text{s}}=\dot{%
\mathbf{R}}_{i}\mathbf{R}_{i}^{T}$ define the body-fixed and spatial angular
velocity. The vector $\mathbf{v}_{i}^{\text{b}}=\mathbf{R}_{i}^{T}\dot{%
\mathbf{r}}_{i}$ is the body-fixed translational velocity, i.e. the velocity
of the origin of $\mathcal{F}_{i}$ measured in the IFR $\mathcal{F}_{0}$ and
resolved in $\mathcal{F}_{i}$. The spatial translational velocity, $\mathbf{v%
}_{i}^{\text{s}}=\dot{\mathbf{r}}_{i}+\mathbf{r}_{i}\times \mathbf{\omega }%
_{i}^{\text{s}}$, is the velocity of the point of the body that is
momentarily passing through the origin of the IFR $\mathcal{F}_{0}$ resolved
in the IFR.

Twists represented in different frames transform as screws according to (\ref%
{Ad}). Let ${^{1}\mathbf{V}_{i}^{1}}=\left( {^{1}}\mathbf{\omega }{_{i}},{%
^{1}}\mathbf{v}{_{i}^{1}}\right) ^{T}$ be the twist of body $i$ represented
in a frame $\mathcal{F}_{1}$. Let $S_{21}=\left( \mathbf{R}_{21},{^{2}}%
\mathbf{r}{_{21}}\right) $ be the frame transformation from $\mathcal{F}_{1}$
to another frame $\mathcal{F}_{2}$. Then ${^{2}\mathbf{V}_{i}^{2}}=\mathbf{Ad%
}_{\mathbf{S}_{21}}{^{1}\mathbf{V}_{i}^{1}}$ is the twist of body $i$
represented in $\mathcal{F}_{2}$. This is ${^{2}\mathbf{V}_{i}^{2}}=\left( {%
^{2}}\mathbf{\omega }{_{i}},{^{2}}\mathbf{v}{_{i}^{2}}\right) ^{T}$, where ${%
^{2}}\mathbf{v}{_{i}^{2}}={^{2}}\mathbf{v}{_{i}^{1}}+{^{2}}\mathbf{r}{_{21}}%
\times {^{2}}\mathbf{\omega }{_{i}}$ is the translational velocity of the
point of body $i$ traveling through the origin of $\mathcal{F}_{2}$, with ${%
^{2}}\mathbf{v}{_{i}^{1}}=\mathbf{R}_{21}{^{1}}\mathbf{v}{_{i}^{1}}$ and ${%
^{2}}\mathbf{\omega }{_{i}}=\mathbf{R}_{21}{^{1}}\mathbf{\omega }{_{i}}$.
The twist is represented, i.e. measured and resolved, in $\mathcal{F}_{2}$.
The vector components can be resolved in yet another frame $\mathcal{F}_{3}$%
. If $\mathbf{R}_{32}$ is the rotation matrix of this change of coordinates,
then ${^{3}\mathbf{V}_{i}^{2}}=\mathbf{Ad}_{\mathbf{R}_{32}}{^{2}\mathbf{V}%
_{i}^{2}}$, with (\ref{AdRr}), is the twist measured in $\mathcal{F}_{2}$
but resolved in $\mathcal{F}_{3}$: ${^{3}\mathbf{V}_{i}^{2}}=\left( {^{3}}%
\mathbf{\omega }{_{i}},{^{3}}\mathbf{v}{_{i}^{2}}\right) ^{T}$, with ${^{3}}%
\mathbf{\omega }{_{i}}=\mathbf{R}_{32}{^{2}}\mathbf{\omega }{_{i},^{3}}%
\mathbf{v}{_{i}^{2}}=\mathbf{R}_{32}{^{2}}\mathbf{v}{_{i}^{2}}={^{3}}\mathbf{%
v}{_{i}^{1}}+{^{3}}\mathbf{r}{_{21}}\times {^{3}}\mathbf{\omega }{_{i}}$.
This is indeed not a frame transformation but rather a coordinate
transformation. Only if twists are measured and resolved in the same frame,
like body-fixed and spatial twists, the screw coordinate transformation is a
frame transformation. In particular, $\mathbf{V}_{i}^{\text{s}}=\mathbf{Ad}_{%
\mathbf{C}_{i}}\mathbf{V}_{i}^{\text{b}}$.

The relative twist of body $j$ w.r.t. body $i$, i.e. the twist of $\mathcal{F%
}_{i}$ represented in $\mathcal{F}_{j}$, is readily defined as ${^{j}%
\widehat{{\mathbf{V}}}_{i}}=\dot{\mathbf{C}}_{j,i}\mathbf{C}_{j,i}^{-1}$.
With (\ref{Cij}) and (\ref{Ad1}) this is%
\begin{eqnarray}
{^{j}\widehat{{\mathbf{V}}}_{i}} &=&\frac{d}{dt}\left( \mathbf{C}_{j}^{-1}%
\mathbf{C}_{i}\right) \mathbf{C}_{i}^{-1}\mathbf{C}_{j}=\mathbf{C}_{j}^{-1}%
\dot{\mathbf{C}}_{i}\mathbf{C}_{i}^{-1}\mathbf{C}_{j}-\mathbf{C}_{j}^{-1}%
\dot{\mathbf{C}}_{j}\mathbf{C}_{j}^{-1}\mathbf{C}_{j}  \notag \\
&=&\mathrm{Ad}_{\mathbf{C}_{j}}^{-1}({\widehat{{\mathbf{V}}}_{i}^{\text{s}}}-%
{\widehat{{\mathbf{V}}}_{j}^{\text{s}})}=\mathrm{Ad}_{\mathbf{C}_{j}}^{-1}({%
\widehat{{\mathbf{V}}}_{i}^{\text{s}})-\widehat{{\mathbf{V}}}_{j}^{\text{b}}}%
.
\end{eqnarray}%
This is the difference of the twists of the two bodies represented in the
BFR at body $j$.

There are yet two further forms of twist used in MBS kinematics. In the 
\emph{hybrid form}, denoted with $\mathbf{V}_{i}^{\text{h}}=\left( \mathbf{%
\omega }_{i}^{\text{s}},\dot{\mathbf{r}}_{i}\right) ^{T}={^{0}{\mathbf{V}}%
_{i}^{i}}$, the twist of body $i$ is measured in the body-fixed frame $%
\mathcal{F}_{i}$ but resolved in the IFR $\mathcal{F}_{0}$. The \emph{mixed
form} of twists, denoted with $\mathbf{V}_{i}^{\text{m}}=\left( \mathbf{%
\omega }_{i}^{\text{b}},\dot{\mathbf{r}}_{i}\right) ^{T}$, uses the
body-fixed angular velocity $\mathbf{\omega }_{i}^{\text{b}}$ and the
translational velocity $\dot{\mathbf{r}}_{i}$. The two forms are related by%
\begin{equation}
\mathbf{V}_{i}^{\text{m}}=\left( 
\begin{array}{cc}
\mathbf{R}_{i}^{T} & \mathbf{0} \\ 
\mathbf{0} & \mathbf{I}%
\end{array}%
\right) \mathbf{V}_{i}^{\text{h}}=\left( 
\begin{array}{cc}
\mathbf{I} & \mathbf{0} \\ 
\mathbf{0} & \mathbf{R}_{i}%
\end{array}%
\right) \mathbf{V}_{i}^{\text{b}}=\left( 
\begin{array}{cc}
\mathbf{R}_{i}^{T} & \mathbf{0} \\ 
-\widetilde{\mathbf{r}}_{i} & \mathbf{I}%
\end{array}%
\right) \mathbf{V}_{i}^{\text{s}}.
\end{equation}%
These transformations are apparently not frame transformations (that would
be described by the adjoint mapping). The definition of twists are
summarized in table \ref{tabTwists}.

\begin{table}[h] \centering%
\begin{tabular}{l|lcc|}
\cline{2-3}\cline{2-4}
& 
\begin{tabular}{l}
Reference point \\ 
is origin of%
\end{tabular}
& 
\begin{tabular}{l}
Frame to resolve \\ 
angular velocity%
\end{tabular}
& 
\begin{tabular}{l}
Frame to resolve \\ 
translational velocity%
\end{tabular}
\\ \hline
\multicolumn{1}{|l|}{$\mathbf{V}_{i}^{\text{s}}$} & \multicolumn{1}{|c}{IFR}
& IFR & IFR \\ 
\multicolumn{1}{|l|}{$\mathbf{V}_{i}^{\text{b}}$} & \multicolumn{1}{|c}{BFR}
& BFR & BFR \\ 
\multicolumn{1}{|l|}{$\mathbf{V}_{i}^{\text{h}}$} & \multicolumn{1}{|c}{BFR}
& IFR & IFR \\ 
\multicolumn{1}{|l|}{$\mathbf{V}_{i}^{\text{m}}$} & \multicolumn{1}{|c}{BFR}
& BFR & IFR \\ \hline
\end{tabular}%
\caption{Summary of the reference point at which the translational velocity
is measured, and the frame in which the angular respectively translational
velocity is resolved}\label{tabTwists}%
\end{table}%

\begin{remark}
It is should be remarked that also the hybrid and mixed twists can be
derived analytically by left and right trivialization, as in (\ref{Vbs}), if
rigid body motions are not considered to be elements of $SE\left( 3\right) $%
. To this end, the direct product group $SO\left( 3\right) \times {\mathbb{R}%
}^{3}$ is used as configuration space of a rigid body: $C=\left( \mathbf{R},%
\mathbf{r}\right) \in SO\left( 3\right) \times {\mathbb{R}}^{3}$. The
multiplication in this group is $C_{1}\cdot C_{2}=(\mathbf{R}_{1}\mathbf{R}%
_{2},\mathbf{r}_{1}+\mathbf{r}_{2})$. This is clearly not a frame
transformation. The hybrid and mixed twist are then $\mathbf{V}^{\text{m}%
}=C^{-1}\cdot \dot{C}=(\mathbf{\omega }^{\text{b}},\dot{\mathbf{r}})$ and $%
\mathbf{V}^{\text{h}}=\dot{C}\cdot C^{-1}=(\mathbf{\omega }^{\text{s}},\dot{%
\mathbf{r}})\in so\left( 3\right) \times \mathbb{R}^{3}$. It must be
emphasized that the \emph{direct product group does not represent screw
motions}. Even though it is occasionally used to model MBS and also in Lie
group integration schemes.
\end{remark}

\newpage%

\section*{B Nomenclature}

\paragraph{\textbf{Symbols:}}

\ 

\begin{tabular}{lllll}
$\mathcal{F}_{0}$ & - Inertial reference frame (IFR) &  &  &  \\ 
$\mathcal{F}_{i}$ & - Body-fixed reference frame (BFR) of body $i$ &  &  & 
\\ 
$\mathcal{J}_{i,i}$ & - Body-fixed joint frame (JFR) for joint $i$ at body $%
i $ &  &  &  \\ 
& \ \ joint $i$ connects body $i$ with its predecessor body $i-1$ &  &  & 
\\ 
$\mathcal{J}_{i-1,i}$ & - JFR for joint $i$ at body $i-1$ &  &  &  \\ 
${^{i}}\mathbf{r}$ & - Coordinate representation of a vector resolved in BFR
on body $i$. &  &  &  \\ 
& \ \ The index is omitted if this is the IFR: $\mathbf{r}\equiv {^{0}}%
\mathbf{r}$. &  &  &  \\ 
$\mathbf{R}_{i}$ & - Rotation matrix from BFR $\mathcal{F}_{i}$ at body $i$
to IFR $\mathcal{F}_{0}$ &  &  &  \\ 
$\mathbf{R}_{i,j}$ & - Rotation matrix transforming coordinates resolved in
BFR $\mathcal{F}_{j}$ &  &  &  \\ 
& \ \ to coordinates resolved in $\mathcal{F}_{i}$ &  &  &  \\ 
$\mathbf{r}_{i}$ & - Position vector of origin of BFR $\mathcal{F}_{i}$ at
body $i$ resolved in IFR $\mathcal{F}_{0}$ &  &  &  \\ 
$\mathbf{r}_{i,j}$ & - Position vector from origin of BFR $\mathcal{F}_{i}$
to origin of BFR $\mathcal{F}_{j}$ &  &  &  \\ 
$\widetilde{\mathbf{x}}$ & - skew symmetric matrix associated to the vector $%
\mathbf{x}\in {\mathbb{R}}^{3}$ &  &  &  \\ 
$C_{i}=\left( \mathbf{R}_{i},\mathbf{r}_{i}\right) 
\hspace{-2ex}%
$ & - Absolute configuration of body $i$. This is denoted in matrix form
with $\mathbf{C}_{i}$ &  &  &  \\ 
$\mathbf{C}_{i,j}=\mathbf{C}_{i}^{-1}\mathbf{C}_{j}%
\hspace{-2ex}%
$ & - Relative configuration of body $j$ w.r.t. body $i$ &  &  &  \\ 
${^{k}}\mathbf{v}{_{i}^{j}}$ & - Translational velocity of body $i$ measured
at origin of BFR $\mathcal{F}_{j}$, resolved in BFR $\mathcal{F}_{k}$ &  & 
&  \\ 
$\mathbf{v}_{i}^{\text{b}}\equiv {^{i}}\mathbf{v}{_{i}^{i}}$ & - Body-fixed
representation of the translational velocity of body $i$ &  &  &  \\ 
$\mathbf{v}_{i}^{\text{s}}\equiv {^{0}}\mathbf{v}{_{i}^{0}}$ & - Spatial
representation of the translational velocity of body $i$ &  &  &  \\ 
${^{k}}\mathbf{\omega }{_{i}}$ & - Angular velocity of body $i$ measured and
resolved in BFR $\mathcal{F}_{k}$ &  &  &  \\ 
$\mathbf{\omega }{_{i}^{\text{b}}}\equiv {^{i}}\mathbf{\omega }{_{i}}$ & -
Body-fixed representation of the angular velocity of body $i$ &  &  &  \\ 
$\mathbf{\omega }{_{i}^{\text{s}}}\equiv {^{0}}\mathbf{\omega }{_{i}}$ & -
Spatial representation of the angular velocity of body $i$ &  &  &  \\ 
${^{k}}\mathbf{V}{_{i}^{j}}$ & - Twist of (RFR of) body $i$ measured in $%
\mathcal{F}_{j}$ and resolved in $\mathcal{F}_{k}$ &  &  &  \\ 
$\mathbf{V}_{i}^{\text{b}}\equiv {^{i}}\mathbf{V}{_{i}^{i}}$ & - Body-fixed
representation of the twist of body $i$ &  &  &  \\ 
$\mathbf{V}{_{i}^{\text{s}}}\equiv {^{0}}\mathbf{V}{_{i}^{0}}$ & - Spatial
representation of the twist of body $i$ &  &  &  \\ 
$\mathbf{V}_{i}^{\text{h}}\equiv {^{0}}\mathbf{V}{_{i}^{i}}$ & - Hybrid form
of the twist of body $i$ &  &  &  \\ 
$\mathsf{V}{^{\text{b}}}$ & - Vector of system twists in body-fixed
representation &  &  &  \\ 
$\mathsf{V}{^{\text{s}}}$ & - Vector of system twists in spatial
representation &  &  &  \\ 
$\mathsf{V}{^{\text{h}}}$ & - Vector of system twists in hybrid
representation &  &  &  \\ 
$\mathsf{V}{^{\text{m}}}$ & - Vector of system twists in mixed representation
&  &  &  \\ 
$\mathbf{Ad}_{\mathbf{R}}$ & - Screw transformation associated with $%
C=\left( \mathbf{R},\mathbf{0}\right) $ &  &  &  \\ 
$\mathbf{Ad}_{\mathbf{r}}$ & - Screw transformation associated with $%
C=\left( \mathbf{I},\mathbf{r}\right) $ &  &  &  \\ 
$\mathbf{Ad}_{\mathbf{C}_{i,j}}$ & - Transformation matrix transforming
screw coordinates represented in $\mathcal{F}_{j}$ &  &  &  \\ 
& \ \ to screw coordinates represented in $\mathcal{F}_{i}$ &  &  &  \\ 
$\widehat{\mathbf{X}}\in se\left( 3\right) $ & - $4\times 4$ matrix
associated with the screw coordinate vectors $\mathbf{X}\in {\mathbb{R}}^{6}$
&  &  &  \\ 
$SE\left( 3\right) $ & - Special Euclidean group in three dimensions -- Lie
group of rigid body motions &  &  &  \\ 
$se\left( 3\right) $ & - Lie algebra of $SE\left( 3\right) $ -- algebra of
screws &  &  &  \\ 
$\mathbf{q}\in {\mathbb{V}}^{n}$ & - Joint coordinate vector &  &  &  \\ 
${\mathbb{V}}^{n}$ & - Configuration space &  &  & 
\end{tabular}

\paragraph{\textbf{Joint Screw Coordinates:}}

\ The introduction of joint screw coordinates requires specification of a
frame in which the screw is measured and a frame where the coordinates are
solved. To this end, the following notation for screw coordinates is used:%
\begin{equation*}
{^{k}\mathbf{X}_{i}^{j}},\ \ \ \ \ \ 
\begin{array}{ll}
i & \text{- index of joint} \\ 
j & \text{- index of frame in which the joint screw is measured} \\ 
k & \text{- index of frame in which the coordinates are resolved}%
\end{array}%
\end{equation*}%
A screw is said to be \emph{represented in frame} $j$ if $j=k$. To simplify
the notation the following short hand notation is used: ${^{j}\mathbf{X}_{i}}%
:={^{j}}\mathbf{X}{_{i}^{j}},\ \mathbf{X}{_{i}}:={^{0}}\mathbf{X}{_{i}}{^{0}}
$.

\begin{tabular}{rllll}
& \  &  &  &  \\ 
& \  &  &  &  \\ 
$^{i-1}\mathbf{Z}_{i}$ & - Constant screw coordinate vector of joint $i$
represented in JFR $\mathcal{J}_{i-1,i}$ at body $i-1$, &  &  &  \\ 
& \ \ determined in the reference configuration $\mathbf{q}=\mathbf{0}$ &  & 
&  \\ 
${^{i}}\mathbf{X}_{i}$ & - Constant screw coordinate vector of joint $i$
represented in BFR $\mathcal{F}_{i}$ at body $i$, &  &  &  \\ 
& \ \ determined in the reference configuration $\mathbf{q}=\mathbf{0}$ &  & 
&  \\ 
${^{i-1}}\bar{\mathbf{X}}_{i}$ & - Constant screw coordinate vector of joint 
$i$ represented in BFR $\mathcal{F}_{i-1}$ at body $i-1$, &  &  &  \\ 
& \ \ determined in the reference configuration $\mathbf{q}=\mathbf{0}$ &  & 
&  \\ 
$\mathbf{Y}_{i}$ & - Constant screw coordinate vector of joint $i$
represented in IFR $\mathcal{F}_{0}$, &  &  &  \\ 
& \ \ determined in the reference configuration $\mathbf{q}=\mathbf{0}$ &  & 
&  \\ 
${^{0}}\mathbf{X}_{i}^{i}$ & - Instantaneous screw coordinate vector of
joint $i$ measured in BFR $\mathcal{F}_{i}$ &  &  &  \\ 
& \ \ but resolved in IFR $\mathcal{F}_{0}$ &  &  & 
\end{tabular}

\paragraph{\textbf{Representations of Velocity of Rigid Body }$i$\textbf{:}}

\ 
\vspace{2ex}%

\begin{tabular}{lllll}
\textbf{Body-fixed twist:} & 
\begin{tabular}{ll}
$\mathbf{V}_{i}^{\text{b}}=\left( 
\begin{array}{c}
\mathbf{\omega }_{i}^{\text{b}} \\ 
\mathbf{v}_{i}^{\text{b}}%
\end{array}%
\right) ,$ & with $\widetilde{\mathbf{\omega }}_{i}^{\text{b}}:={^{i}}%
\widetilde{\mathbf{\omega }}{_{i}}=\mathbf{R}_{i}^{T}\dot{\mathbf{R}}_{i},%
\mathbf{v}_{i}^{\text{b}}:={^{i}}\mathbf{v}{_{i}^{i}}=\mathbf{R}_{i}^{T}\dot{%
\mathbf{r}}_{i}$%
\end{tabular}
&  &  &  \\ 
& 
\begin{tabular}{l}
Twist represented in BFR $\mathcal{F}_{i}$, i.e. measured and resolved in $%
\mathcal{F}_{i}$ \\ 
\ 
\end{tabular}
&  &  &  \\ 
\textbf{Spatial twist:} & 
\begin{tabular}{ll}
$\mathbf{V}_{i}^{\text{s}}=\left( 
\begin{array}{c}
\mathbf{\omega }_{i}^{\text{s}} \\ 
\mathbf{v}_{i}^{\text{s}}%
\end{array}%
\right) ,$ & with $\widetilde{\mathbf{\omega }}_{i}^{\text{s}}:={^{0}%
\widetilde{\mathbf{\omega }}{_{i}}}=\dot{\mathbf{R}}_{i}\mathbf{R}_{i}^{T},%
\mathbf{v}_{i}^{\text{s}}:={^{0}}\mathbf{v}{_{i}^{0}}=\dot{\mathbf{r}}_{i}-%
\widetilde{\mathbf{\omega }}{_{i}^{\text{s}}}\mathbf{r}_{i}$%
\end{tabular}
&  &  &  \\ 
& 
\begin{tabular}{l}
Twist represented in IFR $\mathcal{F}_{0}$, i.e. measured and resolved in $%
\mathcal{F}_{0}$ \\ 
\ 
\end{tabular}
&  &  &  \\ 
\textbf{Hybrid twist:} & 
\begin{tabular}{ll}
$\mathbf{V}_{i}^{\text{h}}=\left( 
\begin{array}{c}
\mathbf{\omega }_{i}^{\text{s}} \\ 
\dot{\mathbf{r}}_{i}%
\end{array}%
\right) $ & 
\end{tabular}
&  &  &  \\ 
& 
\begin{tabular}{l}
Twist measured in BFR $\mathcal{F}_{i}$ but resolved in IFR $\mathcal{F}_{0}$
\\ 
\ 
\end{tabular}
&  &  &  \\ 
\textbf{Mixed twist:} & 
\begin{tabular}{ll}
$\mathbf{V}_{i}^{\text{m}}=\left( 
\begin{array}{c}
\mathbf{\omega }_{i}^{\text{b}} \\ 
\dot{\mathbf{r}}_{i}%
\end{array}%
\right) $ & 
\end{tabular}
&  &  &  \\ 
& 
\begin{tabular}{l}
Angular velocity is measured and resolved in BFR $\mathcal{F}_{i}$. \\ 
Translational velocity is measured in BFR $\mathcal{F}_{i}$ but resolved in
IFR $\mathcal{F}_{0}$.\ 
\end{tabular}
&  &  & 
\end{tabular}

\section*{Acknowledgement}

The author acknowledges that this work has been partially supported by the
Austrian COMET-K2 program of the Linz Center of Mechatronics (LCM).

\end{document}